\documentclass[final]{siamltex}

\usepackage{latexsym}
\usepackage{amssymb,amsbsy,amsmath,amsfonts,amssymb,amscd}
\usepackage{mathrsfs}

\usepackage{float}
\usepackage[caption = false]{subfig}
\usepackage[final]{graphicx}

\usepackage{subfig}
\usepackage{graphicx}

\usepackage{color}

\usepackage{cases}

\usepackage{bm}
\newcommand{\rd}{\mathrm{d}}
\newcommand{\half}{\frac{1}{2}}

\usepackage{subeqnarray}
\usepackage{cases}
\usepackage{xcolor}

\title{\Large \bf Towards understanding the boundary propagation speeds in tumor growth models}
\author{Jian-Guo Liu\thanks{Department of Mathematics and Department of Physics, Duke University} \and
  Min Tang\thanks{School of mathematics, Institute of Natural Sciences and MOE-LSC. Shanghai Jiao Tong University}
\and  Li Wang\thanks{School of Mathematics, University of Minnesota} 
\and Zhennan Zhou \thanks{Beijing International Center for Mathematical Research, Peking University}
}

\begin{document}
\maketitle

\begin{abstract}
At the continuous level, we consider two types of  tumor growth models: the cell density model,  which is based on the fluid mechanical construction, is more favorable for scientific interpretation and numerical simulations; and the free boundary model, as the incompressible limit of the former, is more tractable when investigating the boundary propagation. In this work, we aim to investigate the boundary propagation speeds in those models based on asymptotic analysis of the free boundary model and efficient numerical simulations of the cell density model. We derive, for the first time, some analytical solutions for the free boundary model with pressures jumps across the tumor boundary in multi-dimensions with finite tumor sizes. We further show that in the large radius limit, the analytical solutions to the free boundary model in one and multiple spatial dimensions converge to traveling wave solutions. The convergence rate in the propagation speeds are algebraic in multi-dimensions as opposed to the exponential convergence in 1D. We also propose an accurate front capturing numerical scheme for the cell density model, and extensive numerical tests are provided to verify the analytical findings. 

\end{abstract}

\noindent {\bf Key-words:} Tumor growth models, Brinkman model, free boundary model, front capturing scheme.
\\
\noindent {\bf Mathematics Subject Classification} 35K55; 35B25; 76D27; 76M20; 92C50.

\section{Introduction}


 The invading of solid tumors into a host tissue has been one of the most active areas for mathematical modeling. The tumor density can be influenced by concentration of nutrients, cell division, the extracellular matrix as well as other environmental factors. There are numerous models, including individual-based models, fluid mechanical models, free boundary models, for  tumors in different scenarios \cite{BLM,BOBAM, ByDr,DP,Fredman,Pbook,RCM}. The individual-based model is more accurate for small-scale problems while the later two types of models are built from continuum mechanics \cite{ByDr,DP}. One common question is to understand the propagation speed of tumor boundaries \cite{CFA,Greenspan,Fredman}, and it is also one of the most popular research topics for reaction diffusion equations in general \cite{BH}.

Tumor expands with a constant speed 
has been observed and studied in previous literatures \cite{CFA,TVCVDP}, however such a phenomenon could only be observed for large-scale tumors, which leaves a natural open question: when the tumor size is not large enough, how does the tumor boundary propagates with time? More specifically, this consists of two levels of investigation. One is to figure out the dependence of the limiting constant speed on the model parameters, and the other one is to explore the convergence rate of the propagation speed towards the limit. In particular, it involves a subtler question, whether or not the convergence rate depends on the dimension, as is pointed out in \cite{RBH}. 

In this paper, we investigate two types of continuous models:  the cell density model, which is based on a fluid mechanical construction (see e.g. \cite{Pbook, PQV, PV}), and the free boundary model, which describes the geometric motion of solid tumor borders (see \cite{Fredman} and references therein). The cell density models carry the biggest capacity for scientific interpretations. However, due to the nonlinearity and the lack of analytical solutions, it seems impossible to find the analytical formula of the associated boundary propagation speed. Previously, the convergence of the cell density model to its incompressible limit, which is the free boundary model, has been rigorously justified (\cite{PQV, PV}), but despite the vast interests from both the mathematics and the science communities, analyzing the consistency of the propagation speeds from these models remains at the intuitive level.
In this work, we aim to investigate the connections of the propagation speeds, based on numerical implementations of the cell density model, and asymptotic analysis of the free boundary model.

 We specify the cell density model in the following, which is derived mainly from the assumptions that  the expansion of tumor cells is driven by the cell division and the mechanical pressure \cite{Pbook, RBH}. More precisely,   
we consider the following advection-reaction model as in \cite{PQV,PV,TVCVDP}:
\begin{equation}\label{main}
\frac{\partial }{\partial t } \rho - C_S \nabla \cdot (\rho \nabla W) = \Phi (\Sigma, \rho), \quad \bm x \in \mathbb R^n, \quad t \in \mathbb R^+,
\end{equation}
where $\rho(\bm x,t)$ is the density function of tumor cells, $\Sigma (\rho)$ is the elastic pressure, and $\Phi$ is the growth function. The potential $W$ is related to the pressure $\Sigma$ via the Brinkman model
\begin{equation} \label{eq:Brinkman}
- C_z \Delta W +W =\Sigma,  \quad  \nabla W(\bm x, t) \rightarrow \mathbf 0 \,\, ( |\bm x| \rightarrow + \infty).
\end{equation}
One can write \eqref{main} into 
 the following equation  
\[
\frac{\partial }{\partial t } \rho + \nabla \cdot (\rho \bm v) = \Phi (\Sigma, \rho) \,,
\]
where $\bm v= -C_S \nabla W$ is the velocity field field. The velocity field is curl free and the Brinkman model \eqref{eq:Brinkman} rewrites $-C_z \Delta \bm v + \bm v = - C_S \nabla \Sigma$. When $C_z=0$, the Brinkman model recovers the Darcy's law which says cells move in the direction of  the nagative pressure gradient. A lot of works have been dedicated to case with Darcy's law, see \cite{AP,PQV,Vazquez} and the references therein. When $C_z\neq 0$, the dissipation in velocity due to the internal cell friction is analyzed in \cite{PV}, and the authors have pointed out the theory of mixtures which allows for general formalism combining both the Darcy's law and the Brinkman's law. Similar systems with nonlocal cell interactions and cell growth can be found in \cite{AKY,CKY}.

To complete the cell density model, one has to specify the state equation $\Sigma(\rho)$ and the growth function $\Phi(\Sigma,\rho)$. We assume that the tumor cells are modeled as visco-elastic balls and the elastic pressure is an increasing function of the population density. After neglecting cell adhesion and assuming that $\Sigma(\rho)=0$ when cells are not in contact, one possible choice of the state equation writes   
\begin{equation}  \label{constituion0}
\Sigma = 
\begin{cases}
0, & \rho \le 1\,, \\
C_{\nu} \ln \rho, &  \rho \ge 1. 
\end{cases} 
\end{equation}
The biophysical derivation of \eqref{constituion0} can be found in \cite{TVCVDP}.  It is worth mentioning that other forms of state equations, such as $\Sigma(\rho)=\rho^m$, have been proposed and studied as well \cite{PQV,PTV}.
Let $H$ denote the Heaviside function, i.e. $H(v)=0$ for $v<0$ and $H(v)=1$ for $v>0$, the growth term is chosen to be 
\begin{equation}\label{eq:Phi}
\Phi(\rho)=\rho H\big(C_p-\Sigma(\rho)\big).
\end{equation}
This indicates that when the pressure is less than a threshold denoted by $C_p$, i.e. $\Sigma(\rho)<C_p$, the cell density grow exponentially, while the cell division stops when the precess exceeds the threshold, $\Sigma(\rho)>C_p$. Though the state equation in \eqref{constituion0} and the growth function \eqref{eq:Phi} are not yet experimentally verified, they are qualitatively reasonable and allow for analytical formulations of the front speed.

Similar to \cite{PQV,PTV,PV,TVCVDP},
the fluid mechanical model \eqref{main} \eqref{eq:Brinkman} relates to a free boundary model in the incompressible limit ($C_{\nu} \rightarrow \infty$). The derivation of the corresponding free boundary model from \eqref{main} can be seen in a heuristic way as follows.   
Multiplying equation \eqref{main} by $\Sigma' = \frac{C_{\nu}}{\rho}$ in the support of $\Sigma$, we get
\begin{equation} \label{eq:Sigma}
\frac{\partial }{\partial t } \Sigma  -C_S \nabla \Sigma \cdot \nabla W - C_S C_{\nu} \Delta W = C_{\nu} H.
\end{equation}
Formally, sending $C_{\nu} \rightarrow \infty$ yields the relation $-C_S \Delta W =H$ within the support of $\Sigma$.
Thus, in the incompressible limit, we obtain the complementary relation
\begin{equation} \label{limit00}
\Sigma =0 \quad \text{or} \quad -C_S \Delta W =H.
\end{equation}
Formally, one sees from \eqref{eq:Sigma} that if the initial density is compactly supported, then it remains compactly supported with boundary moving with velocity $v=-C_S\nabla W$, and this completes \eqref{limit00}, the free boundary model. In a similar model, such a limit was proved rigorously \cite{PV}.

This free boundary model has been comprehensively studied in \cite{TVCVDP} by explicitly constructing the 1D traveling wave solutions, which implies constant propagation speed of the tumor borders. However, in principle, the 1D traveling wave solution is relevant only when the tumor radius is approaching infinity, and therefore is unable to quantify the dynamics for finite size tumors. It is worth mentioning that the traveling solutions are also available for some multi-species models (see e.g. \cite{LLP}), which sheds light on the understanding of the tumor boundary instability.

In this paper, we construct close-form radially symmetric solutions of the free boundary model in various dimensions. The derivation follows similar techniques to those in \cite{LTWZ19}, but to the best of our knowledge, the results with exact quantification of the pressure jumps are obtained for the first time.  In addition, the expressions provide a strong evidence for the conjecture that the pressure jump relates to the tumor border curvature. We further carry out asymptotic analysis of the close-form solutions in the large tumor radius regime, and are able to identify the traveling wave solutions in the large radius limit. Besides, the asymptotic analysis manifests the effect of the dimension in the large radius limit.
We show that in contrast to the exponential convergence of the speed towards the limit in 1D,  the convergence rate in multi-dimensional cases are at most algebraic.

For the cell density model, direct analysis of the boundary moving speed still seems inaccessible at this stage. Instead, we provide some a priori analysis and propose a novel numerical scheme to simulate its dynamics. The numerical method is an improved version of our previous work \cite{LTWZ18}, wherein only the Darcy's law is considered. When $C_\nu$ is large, though we could only show at a formal level the convergence of the cell density model to the free boundary model, the numerical results show that the boundary moving speed and pressure jump across the tumor borders agree well with the analytical results from the free boundary model. 

The rest of the paper is outlined as follows. We give a priori $L^2$ estimate for the fluid mechanical model in section 2 and then in section 3, based on the limiting free boundary model, we derive explicitly the velocity and structure of the tumor boundary for 1D symmetric, 2D radial symmetric and 3D spherical symmetric cases. A new numerical scheme that captures the correct border velocity for a wide range of $C_{\nu}$ 
is proposed in section 4, and in section 5 we carry out extensive numerical tests to verify the analytical observations.

\section{A priori analysis of the cell density model}

In the section, we aim to derive some a priori estimates of the cell density model. Note that, although quantifying the boundary propagation speed at this level seems unreachable, the stability results we obtained below gives access to the design of reliable numerical schemes. 

For simplicity, we take all the parameters except $C_\nu$ equal to $1$, i.e., $C_z= C_S=C_p=1$, and the model simplifies to
\begin{align}
\frac{\partial }{\partial t } \rho -  \nabla \cdot (\rho \nabla W) & = \rho H(1 - \Sigma),  \label{eq:thos}
\\
- \Delta W +W & =\Sigma. \label{eq:Ws}
\end{align}
From the constitutive law \eqref{constituion0} it is clear that 
\begin{equation} \label{est1}
0 \le \Sigma(\bm x, t) \le 1, \quad \text{if} \quad  0 \le \Sigma(\bm x, 0) \le 1.
\end{equation}
By the maximum principle, we also have 
$0 \le W(\bm x, t) \le 1$, and from \eqref{eq:Ws}, it implies $
-1 \le \Delta W \le 1$. 
Note that, $ \Delta W$ may change signs in the whole space.

Next, we check the $L^2$ stability of the density $\rho$. Assume $\rho(\bm x,0)$ is compactly supported and $\| \rho(\cdot,0)\|_{L^2}$ is finite, multiplying equation \eqref{eq:thos} by $\rho$ and integrating over $\mathbb R^n$, we get
\[
\frac{1}{2} \frac{d}{dt}   \int_{\mathbb R^d} |\rho|^2 d \bm x  =  \int_{\mathbb R^d} \rho \nabla  \cdot (\rho \nabla W) d \bm x + \int_{\mathbb R^d} |\rho|^2 H(1-\Sigma) d \bm x = 
\frac 1 2 \int_{\mathbb R^n} |\rho|^2 \Delta W d \bm x +   \int_{\mathbb R^d} |\rho|^2 H(1-\Sigma) d \bm x.
\]
which implies that $\frac{d}{dt} \| \rho \|^2_{L^2}  \le {3}\| \rho \|^2_{L^2}$, and therefore upper bound for the relative 
growth rate of $\| \rho\|_{L^2}$ is guaranteed. 

We now analyze the $L^2$ stability of the pressure function $\Sigma$. We denote by $D(t) = \{ \bm x: \Sigma (\bm x, t) \ge 0 \}$, and assume that $D(t)$ is compactly supported in $\mathbb R^n$, then 
\[
\frac{\partial }{\partial t } \Sigma  -\nabla \Sigma \cdot \nabla W - C_\nu \Delta W =C_\nu H, \quad x \in D(t).
\]
Multiply each side by $\Sigma$ and integrate over $\mathbb R^n$, we have
\[
\frac{d}{dt} \int_{\mathbb R^n} |\Sigma|^2 d \bm x - \int_{\mathbb R^n} \Sigma \nabla \Sigma \cdot \nabla W  d \bm x - C_\nu\int_{\mathbb R^n} \Sigma \Delta W  d \bm x= C_\nu \int_{D(t)} \Sigma H d \bm x. 
\]
By equation \eqref{eq:Ws} and the boundedness of $W$ and $\Delta W$ above, we get
\begin{align*}
\int_{\mathbb R^n} \Sigma \nabla \Sigma \cdot \nabla W  d \bm x& = -\frac{1}{2} \int_{\mathbb R^n} (\Sigma)^2 \Delta W d \bm x \\
& = - \frac{1}{2} \int_{\mathbb R^n} (- \Delta W +W)^2 \Delta W d \bm x  \\
& =  - \frac{1}{2} \int_{\mathbb R^n} \left[ |\Delta W|^2 \Delta W- 2 |\Delta W|^2 W + W^2 \Delta W  \right] d \bm x \\
& = - \frac{1}{2} \int_{\mathbb R^n} \left[ |\Delta W|^2 \Delta W- 2 |\Delta W|^2 W - 2W | \nabla W|^2  \right] d \bm x  \\
& \le \,\, \frac 3 2 \int_{\mathbb R^n} |\Delta W |^2 d \bm x + \int_{\mathbb R^n} |\nabla W |^2 d \bm x,
\end{align*}
and 
\begin{align*}
 C_\nu\int_{\mathbb R^n} \Sigma \Delta W  d \bm x &=  C_\nu\int_{\mathbb R^n} (- \Delta W +W) \Delta W d \bm x = - C_\nu \int_{\mathbb R^n} |\Delta W|^2 + |\nabla W|^2 d \bm x.
\end{align*}
Denote $V(t)=\text{Vol}(D(t))$, then with \eqref{est1} we obtain $C_\nu \int_{D(t)} \Sigma H d \bm x \le C_\nu V(t)$. Altogether, we have the following estimate
\begin{equation}
\frac{d}{dt} \int_{\mathbb R^n} |\Sigma|^2 d \bm x \le \left(\frac 3 2 - C_\nu \right) \int_{\mathbb R^n} |\Delta W |^2 d \bm x + (1-C_\nu) \int_{\mathbb R^n} |\nabla W |^2 d \bm x+C_\nu V(t).
\end{equation}
Clearly, when $C_\nu > \frac{3}{2}$, diffusion dominates the convection and results in an overall stabilizing effect. 




\section{The free boundary model} \label{sec:fbm}

In this section, we construct analytical solutions to the free boundary problem based on the three-zone ansatz, which was originally proposed in \cite{TVCVDP} for the constructing of traveling wave solutions. However, unlike the traveling wave solution where the inner layer is infinite, here we assume that the inner layer has a finite size. We shall show that, with the specific choice of solution ansatz described below, the free boundary model reduces to a differential-algebraic system of equations, where the differential equation of the radius determines the border expanding speed and the algebraic equation governs the thicknesses of the inner and outer layers of the tumor. 

We will also investigate the large radius limit when the thickness of the inner layer becomes infinity. In such limits, the radial symmetric solutions to the free boundary model always converges to a traveling wave solution, regardless of the spatial dimensions, but with different convergence rate. In multi-dimensions the convergence rate are algebraic with respect to the radius of the inner layer, which gives hints to the dependence of the curvature (the reciprocal of the radius in this case) in the first order correction of the front moving speed.

In the incompressible limit, we consider the Hele-Shaw type complementary equation and the Brinkman model
\begin{subequations}  \label{limit000}
\begin{numcases}{}
\Sigma =0 \quad \text{or} \quad -C_S \Delta W =H.
\\ - C_z \Delta W +W =\Sigma.
\end{numcases}
\end{subequations}
The free boundary model is completed by boundary moving velocity $v=-C_S\nabla W$.
Particularly, we are interested in solution with density $\rho$ evolving as a characteristic function on a changing domain and pressure $\Sigma$ may vary within the support of $\rho$. 

\begin{figure}[h!] 
\centering 
\includegraphics[width=0.55\textwidth]{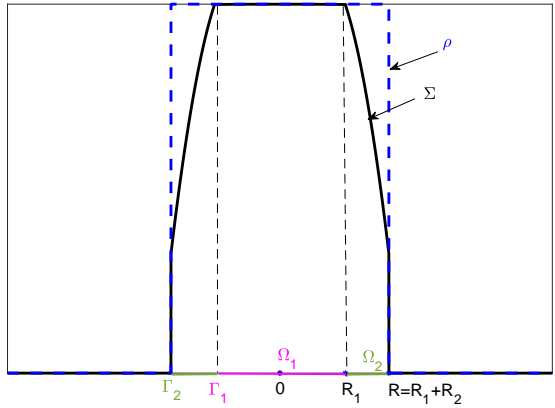}
\caption{Schematic plot of the solution ansatz. The support of the tumor consists of two layers: $\Omega_1$ and $\Omega_2$, which define the inner boundary $\Gamma_1$ and the outer boundary $\Gamma_2$ respectively. Note that, there is a possible pressure jump across $\Gamma_2$.}
\label{fig:3zone}
\end{figure}

To this end, we assume that the whole domain can be divided into three parts: in $\Omega_1$, $\Sigma=C_p$, $\rho=1$, and its boundary is denoted by $\Gamma_1$; in $\Omega_2$, $\Sigma\in (0, C_p)$, $\rho=1$, its inner and outer boundary are $\Gamma_1$ and $\Gamma_2$ respectively, and $\Omega_1 \cap \Omega_2=\emptyset$; $\Omega_3= (\Omega_1 \cup \Omega_2)^c$, where $\Sigma=0$, and $\rho=0$. $\Gamma_2$ is evolving in time with the normal velocity $- \nabla W \cdot \hat n$, where $n$ is the outer unit normal vector to $\Gamma_2$. Please see Figure~\ref{fig:3zone} for illustration.

Note specifically that the support of $\Sigma$ may not coincide with that of $\rho$ in general, but we are only interested in deriving the analytical solutions when they do share the same support. We also expect that, $W$ are $\nabla W$ are continuous across both $\Gamma_1$ and $\Gamma_2$, whereas pressure $\Sigma$ remains continuous across $\Gamma_1$ but has a jump across $\Gamma_2$. We also note that  $W$ is not supported in $\Omega$.  

Since the Heaviside function $H$ is hard to deal with, we adopt the following regularization as in \cite{TVCVDP}: 
\begin{equation} \label{Heta}
H_{\eta}(u)= 
\begin{cases}
0, & u\le0; \\
\frac{u}{\eta}, & 0\le u \le \eta; \\
1, & u \ge \eta \,,
\end{cases}
\end{equation}
where $\eta \in (0, C_p)$. As a result, the decomposition of the domain is modified accordingly. In $\Omega_1^{\eta}$, $\Sigma \in ( C_p -\eta, C_p ]$  and $\rho=1$; in $\Omega_2^{\eta}$, $\Sigma\in (0,C_p-\eta)$ and $\rho=1$; finally in $\Omega_3^{\eta}$, $\Sigma=0$ and $\rho=0$. The continuity of $W$, $\nabla W$ and $\Sigma$ through the boundaries $\Gamma_1^{\eta}$ and $\Gamma_2^{\eta}$ stay  unchanged. It is expected that the regularized solution converges to the original one in the limit $\eta \rightarrow 0^+$.

In the rest of the section, we derive explicit solutions to the incompressible limit model using the above ansatz in various dimensions. We also investigate the solvability conditions in order for the ansatz to be valid, and its connection to the traveling wave solution.


\subsection{1D case}

We start with the regularized problem. For simplicity, we assume the problem is symmetric in space, and denote
$\Omega_1^{\eta}= [-R^{\eta}_1(t), R^{\eta}_1(t)], \quad \Omega_1^{\eta} \cup \Omega_2^{\eta} = [-R^{\eta}(t), R^{\eta}(t)].
$
with
$
R^{\eta} (0)= R_0
$
being the initial condition. We first derive the equation that links $R^{\eta}$ and $R_1^{\eta}$, and then evolution equation for $R^{\eta}$. 

In $\Omega_1^\eta$, \eqref{limit000} along with \eqref{Heta} writes 
\[
-C_S  W_{xx} = \frac{C_p -\Sigma}{\eta}, \qquad 
- C_z W_{xx} +W =\Sigma\,,
\]
which readily leads to
\[
-(\eta C_S+ C_z) W_{xx}+W=C_p\,,
\]
after eliminating $\Sigma$. 
Note from the symmetric assumption that $W'(0)=0$, then the general solution of $W$ in $\Omega_1^\eta$ is given by
\begin{equation*}
W(x)=C_p + A^\eta \cosh (\frac{x}{\sqrt{\eta C_S+C_z}}), \quad x \in \Omega_1^\eta\,.
\end{equation*} 
Consequently, the general solution of $\Sigma$ in $\Omega_1^\eta$ is given by
\begin{equation*}
\Sigma=- C_z W_{xx} +W=C_p + \frac{A^\eta \eta C_S}{\eta C_S+ C_z} \cosh \left(\frac{x}{\sqrt{\eta C_S+C_z}}\right), \quad x \in \Omega_1^\eta. 
\end{equation*}
Since $\Sigma$ equals $C_p-\eta$ at the boundary of $\Omega_1^\eta$, i.e., $\Sigma(R_1^{\eta})= C_p- \eta$, we have
which leads to 
\begin{equation} \label{A-coeff}
A^\eta= - \frac{\eta C_S+C_z}{C_S \cosh \left(\frac{R_1^{\eta}}{\sqrt{\eta C_S+C_z}}\right) }.
\end{equation}

In $\Omega_2^\eta$, the model \eqref{limit000} becomes 
\[
-C_S  W_{xx} = 1, \qquad 
- C_z W_{xx} +W =\Sigma\,,
\]
and one can immediately write down the the general solution for $W$ as 
\begin{equation*}
W(x) = - \frac{1}{2C_S} x^2 + a^\eta x +b^\eta, \qquad x \in \Omega_2^\eta\,.
\end{equation*}
By continuity of $W$ and $W_x$ at $x=R_1^{\eta}$, we get
\begin{equation}\label{a-coeff}
a^\eta=\frac{R_1^{\eta}}{C_S}- \frac{\sqrt{\eta C_S+C_z}}{C_S} \tanh \left(\frac{R_1^{\eta}}{\sqrt{\eta C_S+C_z}}\right), 
\end{equation}
\begin{equation}\label{b-coeff}
b^\eta=C_p - \eta - \frac{C_z}{C_S}- \frac{(R^{\eta}_1)^2}{2C_S} +R_1^{\eta} \frac{\sqrt{\eta C_S+C_z}}{C_S} \tanh \left(\frac{R_1^{\eta}}{\sqrt{\eta C_S+C_z}}\right).
\end{equation}
And the general solution of $\Sigma$ in $\Omega_2^\eta$ is given by
\begin{equation*}
\Sigma(x) =- C_z W_{xx} +W= - \frac{1}{2C_S} x^2 + a^\eta x +b^\eta + \frac{C_z}{C_S}.
\end{equation*}

Finally, in $\Omega_3^\eta$, \eqref{limit000} simplifies to
\[
\Sigma=0, \qquad 
- C_z W_{xx} +W =\Sigma.
\]
By assuming the decaying behavior of $W$ at infinity, the general solution of $W$ in $\Omega_3^\eta$ can be written as 
\begin{equation*}
W(x) =d^\eta e^{- \frac{x-R^{\eta}}{\sqrt{C_z}}}, \qquad x \in \Omega_3^h\,.
\end{equation*}
Then the continuity of $W$ at $R^\eta$ implies
\begin{equation} \label{d-coeff}
d^\eta=- \frac{1}{2C_S} (R^\eta)^2 + a^\eta R^{\eta} +b^\eta.
\end{equation}
To summarize, we have the following analytical representation of $\Sigma$ and $W$ in different domains 
\begin{equation} \label{W-1D}
W(x) = \left\{ \begin{array}{cc} C_p + A^\eta \cosh (\frac{x}{\sqrt{\eta C_S+C_z}}), & x \in \Omega_1^\eta\,,
\\ - \frac{1}{2C_S} x^2 + a^\eta x +b^\eta, & x \in \Omega_2^\eta\,,
\\ d^\eta e^{- \frac{x-R^{\eta}}{\sqrt{C_z}}}\,, & x \in \Omega_3^\eta\,. \end{array}
\right.
\end{equation}
\begin{equation}\label{Sigma-1D}
\Sigma(x) = \left\{ \begin{array}{cc} C_p + \frac{A^\eta \eta C_S}{\eta C_S+ C_z} \cosh \left(\frac{x}{\sqrt{\eta C_S+C_z}}\right), & x \in \Omega_1^\eta\,,
\\ - \frac{1}{2C_S} x^2 + a^\eta x +b^\eta + \frac{C_z}{C_S}, & x \in \Omega_2^\eta\,,
\\ 0 & x \in \Omega_3^\eta\,, 
\end{array}
\right.
\end{equation}
where $a^\eta$, $b^\eta$, $d^\eta$ and $A^\eta$ are given by \eqref{a-coeff}, \eqref{b-coeff}, \eqref{d-coeff} and \eqref{A-coeff} respectively. 

Thus, the regularized problem has been completed solved. We take the limit $\eta \rightarrow 0$, and the solution becomes 
\begin{equation} \label{W-1D-lim}
W(x) = \left\{ \begin{array}{cc} C_p + A \cosh (\frac{x}{\sqrt{C_z}}), & x \in \Omega_1\,,
\\ - \frac{1}{2C_S} x^2 + a x +b, & x \in \Omega_2\,,
\\ de^{- \frac{x-R}{\sqrt{C_z}}}\,, & x \in \Omega_3\,; \end{array}
\right.
\end{equation}
\begin{equation}\label{Sigma-1D-lim}
\Sigma(x) = \left\{ \begin{array}{cc} C_p, & x \in \Omega_1\,,
\\ - \frac{1}{2C_S} x^2 + a x +b + \frac{C_z}{C_S}, & x \in \Omega_2 \,,
\\ 0 & x \in \Omega_3\,, 
\end{array}
\right.
\end{equation}
where the parameters are listed below
\begin{align*}
&a=\frac{R_1}{C_S}- \frac{\sqrt{C_z}}{C_S} \tanh \left(\frac{R_1}{\sqrt{C_z}}\right), \quad
b=C_p - \eta - \frac{C_z}{C_S}- \frac{(R_1)^2}{2C_S} +R_1 \frac{\sqrt{C_z}}{C_S} \tanh \left(\frac{R_1}{\sqrt{C_z}}\right)\,,
\\& A= - \frac{C_z}{C_S \cosh \left(\frac{R_1}{\sqrt{C_z}}\right) }, \qquad 
d=- \frac{1}{2C_S} (R)^2 + a R +b\,.
\end{align*}

\begin{figure}[h!] 
\centering
\includegraphics[width=0.32\textwidth]{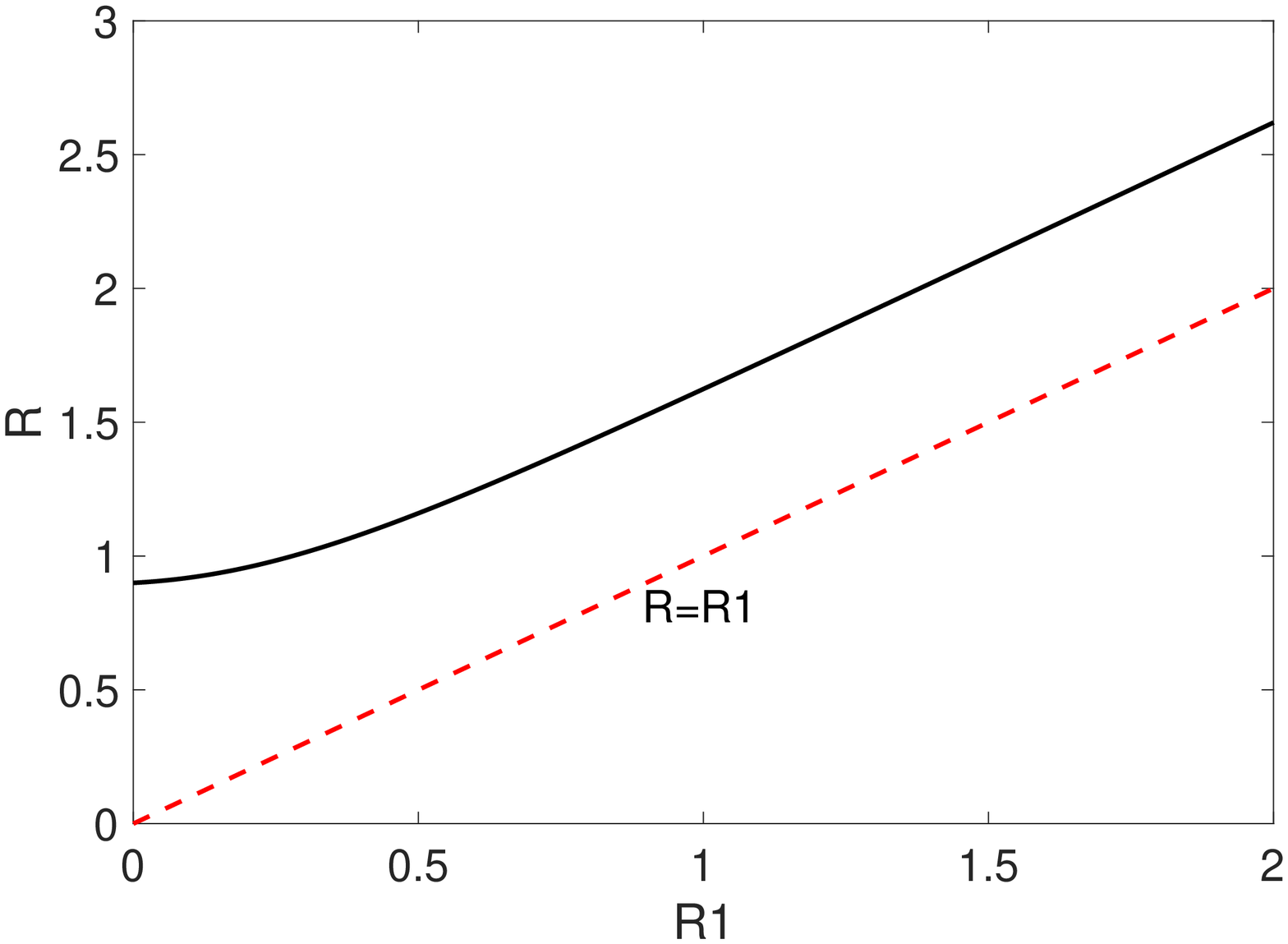}
\includegraphics[width=0.31\textwidth]{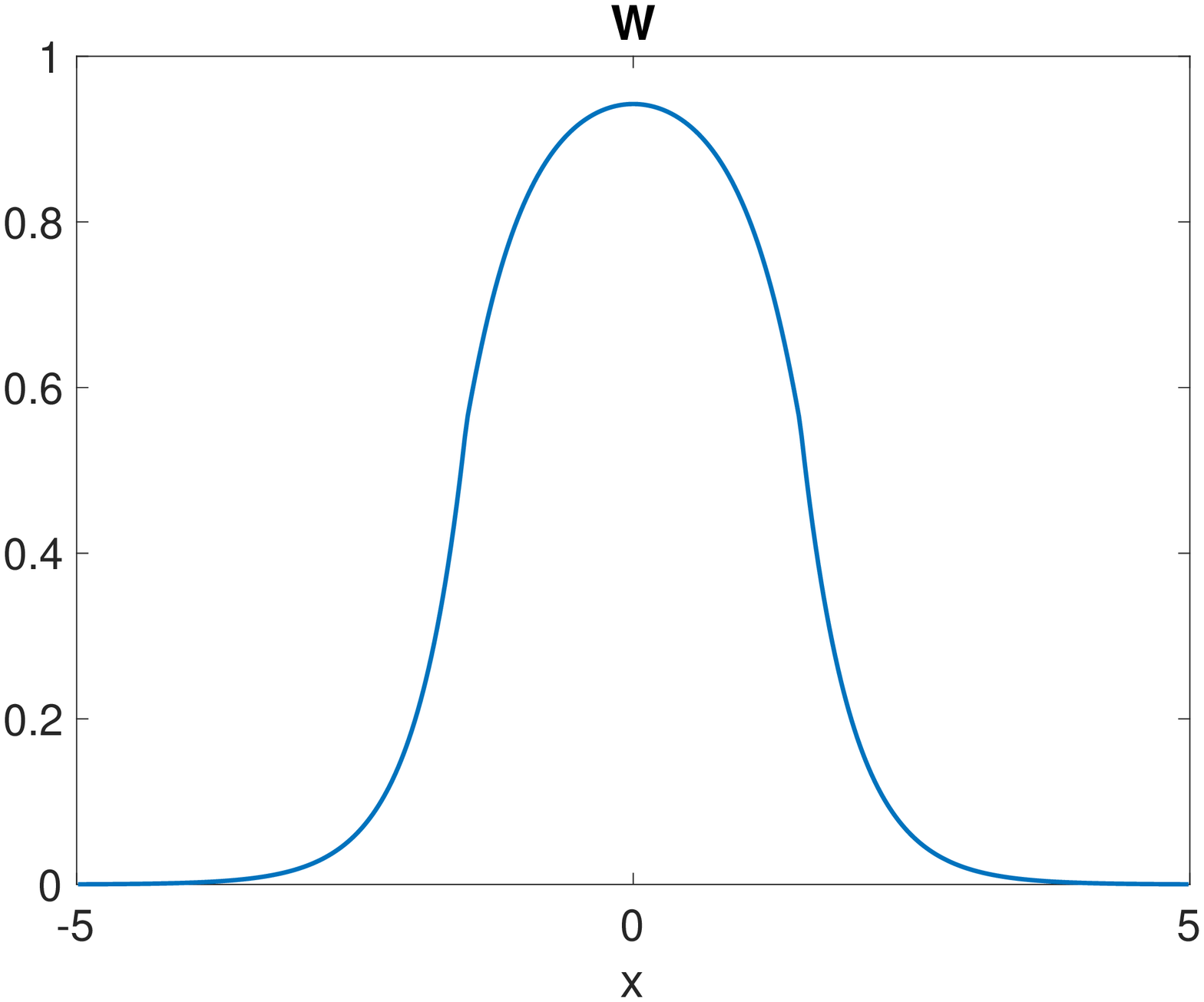}
\includegraphics[width=0.32\textwidth]{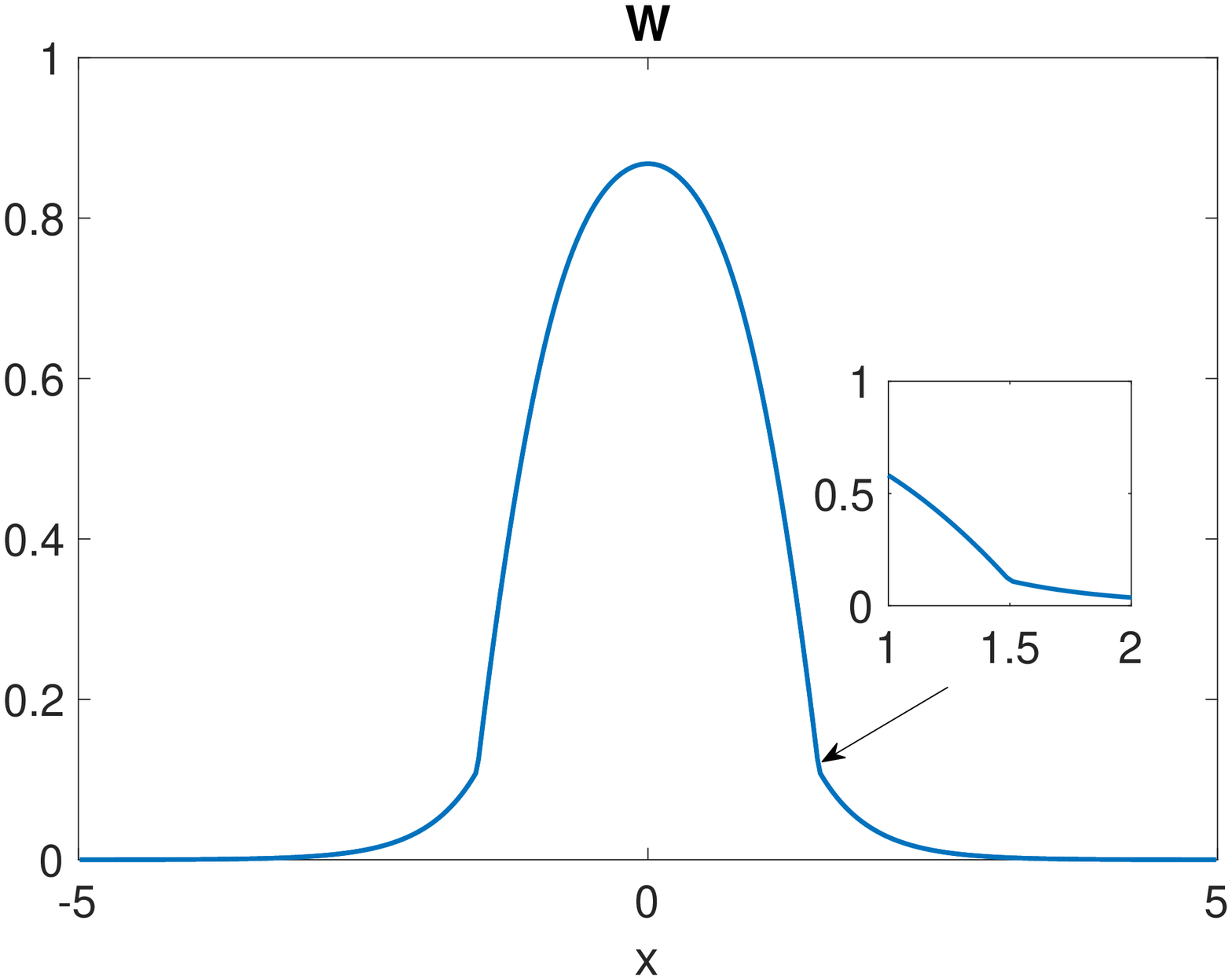}
\caption{Left: relationship between $R_1$ and $R$. Middle: plot of $W(x)$ using \eqref{W-1D-lim} with $R = 1.5$ and $R_1 = 1.2781$, satisfying the relation \eqref{eq:R2}. Right: plot of $W(x)$ using \eqref{W-1D-lim} with $R = 1.5$ and $R_1 = 0.8$, not satisfying the relation \eqref{eq:R2}. 
Here $C_S = 1$,  $C_z = 2$,  $C_p = 4$,  $C_{\nu}= 200$, $\eta = 0.001$.}
\label{R-R2-1D}
\end{figure}

Next, we examine the relationship between two boundaries $R$ and $R_1$. Again by continuity of $W_x$ at $R$, one has 
\begin{equation}\label{rel:R1R}
- \frac{R}{C_S}+ a = - \frac{d}{\sqrt{C_z}}. 
\end{equation}
If we denote the difference between those two by $R_2$, namely $R=R_1+R_2$, then \eqref{rel:R1R} becomes
\[
\frac{\sqrt{C_z}}{C_S}R_2  + \frac{C_z}{C_S} \tanh\left( \frac{R_1}{\sqrt{C_z}} \right) 
= - \frac{1}{2C_S}(R_2)^2 - R_2 \frac{\sqrt{C_z}}{C_S} \tanh\left( \frac{R_1}{\sqrt{C_z}} \right)+ C_p- \frac{C_z}{C_S},
\]
which simplifies to a quadratic equation in $R_2$
\begin{equation} \label{eq:R2}
(R_2)^2+2 \sqrt{C_z} \left( 1+  \tanh\left( \frac{R_1}{\sqrt{C_z}} \right)  \right) R_2+2 C_z \left( 1+  \tanh\left( \frac{R_1}{\sqrt{C_z}} \right)  \right) -2C_pC_S =0.
\end{equation}
We remark that, given the parameters $C_p$, $C_z$ and $C_S$ (note that $C_\nu \rightarrow \infty$), the necessary condition for above solution to make sense is that $R \ge R_1 \ge 0$. This implies that, in the algebraic equation \eqref{eq:R2}, given the value for the outer boundary, $R=R_1+R_2$, there exist solutions with $R_1 \ge 0$ and $R_2 \ge 0$. Although one cannot get explicit constraints from such solvability conditions, it is easy check the condition numerically, see Figure \ref{R-R2-1D} on the left. The rest two plots in Figure \ref{R-R2-1D} indicates that if we choose $R$ and $R_1$ satisfying the relation \eqref{eq:R2}, then $W$ has a smooth transition from $\Omega_1$ to $\Omega_2$ (middle plot), and a kink otherwise (right plot). 

To make a connection to the traveling wave solutions, we consider the limit $R_1 \rightarrow \infty$. Then \eqref{eq:R2} reduces to 
\begin{equation*} 
(R_2)^2+4 \sqrt{C_z}  R_2+4 C_z -2C_pC_S =0,
\end{equation*}
which has the following two solutions
$
R_2= \pm \sqrt{2C_pC_S} - 2\sqrt{C_z}.
$
Therefore, when $\sqrt{2C_pC_S} > 2\sqrt{C_z}$, $R_2$ has a positive solution, which means in the traveling wave limit, i.e. $R_1 \rightarrow \infty$, $\Omega_2$ persists with width $\sqrt{2C_pC_S} - 2\sqrt{C_z}$.

In this case, we can further calculate the pressure jump at $R$, which is given by $\Sigma (R)= d + \frac{C_z}{C_S}$. With \eqref{rel:R1R}, it becomes
\[
\Sigma (R)= \frac{\sqrt{C_z}}{C_S}R_2+  \frac{\sqrt{C_z}\sqrt{C_z}}{C_S} \tanh\left( \frac{R_1}{\sqrt{C_z}} \right)+ \frac{C_z}{C_S}.
\]

If we take $R_1 \rightarrow \infty$, the pressure at $\Gamma_2$  simplifies to 
\[
\Sigma (R)= \frac{\sqrt{C_z}}{C_S}R_2 +2 \frac{C_z}{C_S}.
\]
When $\sqrt{2C_pC_S} > 2\sqrt{C_z}$, we substitute $R_2$ with $\sqrt{2C_pC_S} - 2\sqrt{C_z}$, and get the pressure jump
\[
\Sigma (R)= \sqrt{\frac{2C_zC_p}{C_S}}.
\]
This recovers traveling wave solutions found in \cite{TVCVDP}.


When $R_1$ is finite, there is no explicit analytical solution, instead we numerically solve \eqref{eq:R2} for $R_1$. Here $R(t)$ satisfies
\[
\dot R = - C_S W_x(R) = R-C_S \left( \frac{R_1}{C_S}- \frac{\sqrt{C_z}}{C_S} \tanh \left(\frac{R_1}{\sqrt{C_z}}\right) \right) = R-R_1 + \sqrt{C_z}  \tanh \left(\frac{R_1}{\sqrt{C_z}}\right).
\]
When $\sqrt{2C_pC_S} > 2\sqrt{C_z}$, in the traveling wave limit, $R_1 \rightarrow \infty$, we obtain
\[
\dot R = \sqrt{2C_pC_S} - \sqrt{C_z}.
\]

\subsection{2D radial symmetric case}

Similar to the 1D case, with the radial symmetric assumption, we can explicitly solve for the ansatz solution to the regularized impressible model. The interested readers may refer to Appendix \ref{2d:cal} for details. And by taking the limit $\eta \rightarrow 0$, we obtain the following solution to the incompressible limit model
\begin{equation}\label{W-R-2d-lim}
W(r) = \left\{ \begin{array}{cc} 
C_p + A  I_0 \left(\frac{r}{\sqrt{C_z}}\right), & r\in \Omega_1
\\ - \frac{1}{4C_S} r^2 + a \ln r +b, &  r \in \Omega_2
\\ d  K_0\left(\frac{r}{\sqrt{C_z}} \right), & r \in \Omega_3
\end{array} \right.
\end{equation}
\begin{equation}\label{Sigma-R-2d-lim}
\Sigma(r) = \left\{ \begin{array}{cc} 
C_p , & r\in \Omega_1
\\ - \frac{1}{4C_S} r^2 + a \ln r +b + \frac{C_z}{C_S}, &  r \in \Omega_2
\\ 0, & r \in \Omega_3\,,
\end{array} \right.
\end{equation}
where the parameters are listed below
\begin{align*} 
&a= \frac{1}{2C_S}R_1^2 - R_1 \frac{\sqrt{C_z}}{C_S} \frac{ I_1 \left(\frac{R_1}{\sqrt{C_z}}\right)}{ I_0 \left(\frac{R_1}{\sqrt{C_z}}\right)},
\qquad 
b=C_p - \frac{C_z}{C_S}+ \frac{(R_1)^2}{4C_S} -\frac{(R_1)^2 \ln R_1}{2C_S}  +R_1 \ln R_1 \frac{\sqrt{C_z}}{C_S} \frac{ I_1 \left(\frac{R_1}{\sqrt{C_z}}\right)}{ I_0 \left(\frac{R_1}{\sqrt{C_z}}\right)},
\\ &A = - \frac{C_z}{C_S I_0 \left(\frac{R_1}{\sqrt{C_z}}\right) }, 
\qquad d K_0 \left(\frac{R}{\sqrt{C_z}} \right) =- \frac{1}{4C_S} R^2 + a \ln R +b\,.
\end{align*}

Next, we examine the relationship between two boundaries $R$ and $R_1$. Again by continuity of $W_r$ at $R$, one has
\begin{equation*}
 - \frac{d}{\sqrt{C_z}}K_1 \left(\frac{R}{\sqrt{C_z}} \right) = - \frac{1}{2C_S} R + \frac{a}{R}.
\end{equation*}
If we denote the difference between those two by $R_2$, namely $R=R_1+R_2$, then 
\begin{equation}\label{rel:2d3}
\sqrt{C_z} \frac{K_0\left(\frac{R}{\sqrt{C_z}} \right)}{K_1\left(\frac{R}{\sqrt{C_z}} \right)} \left( \frac{1}{2C_S} R- \frac{a}{R} \right) =- \frac{1}{4C_S} R^2 + a \ln R +b.
\end{equation}

To see the connection to the traveling wave model, we consider the case with $R_1 \gg 1, R_2=O(1)$.
Asymptotically expanding \eqref{rel:2d3}, we get 
\begin{align*}
& \, \, \text{L.H.S.} \\ 
=  &  \sqrt{C_z}  \left(1 - \frac{\sqrt{C_z}}{2R}  \right)\left( \frac{1}{C_S}R_2 + \frac{\sqrt{C_z}}{C_S} \frac{I_1 \left( \frac{R_1}{\sqrt{C_z}} \right)}{I_0 \left( \frac{R_1}{\sqrt{C_z}} \right)} + \frac{1}{R_1} \left( \frac{R_2^2}{2C_S}- R_2 \frac{\sqrt{C_z}}{C_S} \frac{I_1 \left( \frac{R_1}{\sqrt{C_z}} \right)}{I_0 \left( \frac{R_1}{\sqrt{C_z}} \right)} \right)  \right) + o \left( \frac{1}{R_1} \right) \\
= & \sqrt{C_z}  \left(1 - \frac{\sqrt{C_z}}{2R_1}  \right) \left( \frac{1}{C_S}R_2 + \frac{\sqrt{C_z}}{C_S}  + \frac{1}{R_1} \left( \frac{R_2^2}{2C_S}- R_2 \frac{\sqrt{C_z}}{C_S}  - \frac{C_z}{2C_S} \right) \right)+ o \left( \frac{1}{R_1} \right) .
\end{align*}
Here, we have used the fact that, when $z \gg 1$
\[
\frac{I_1(z)}{I_0(z)}= 1 - \frac{1}{2z}+ o\left( \frac 1 z \right), \quad \frac{K_0(z)}{K_1(z)}= 1 - \frac{1}{2z}+ o\left( \frac 1 z \right).
\]
Similarly, on the right hand side, we have
\begin{align*}
& \, \, \text{R.H.S.} \\
= & - \frac{1}{2C_S}R_2^2 - R_2 \frac{\sqrt{C_z}}{C_S} \frac{I_1 \left( \frac{R_1}{\sqrt{C_z}} \right)}{I_0 \left( \frac{R_1}{\sqrt{C_z}} \right)} + C_p - \frac{C_z}{C_S} + \frac{1}{R_1} \left( \frac{R_2^3}{6C_S} + \frac{R_2^2}{2}\frac{\sqrt{C_z}}{C_S} \frac{I_1 \left( \frac{R_1}{\sqrt{C_z}} \right)}{I_0 \left( \frac{R_1}{\sqrt{C_z}} \right)}\right)+ o \left( \frac 1 {R_1} \right) \\
= &  - \frac{1}{2C_S}R_2^2 - R_2 \frac{\sqrt{C_z}}{C_S}  + C_p - \frac{C_z}{C_S} + \frac{1}{R_1} \left( \frac{R_2^3}{6C_S} + \frac{R_2^2}{2}\frac{\sqrt{C_z}}{C_S}+ \frac{R_2 C_z}{2C_S} \right)+ o \left( \frac 1 {R_1} \right).
\end{align*}

To match the terms by order, we assume that when $R_1 \gg 1$, 
\[
R_2 = \alpha_0 + \frac{\alpha_1}{R_1} + o\left( \frac 1 {R_1} \right).
\] Then to the leading order, we have
\[
(\alpha_0)^2+ 4 \sqrt{C_z}  \alpha_0 +4C_z-2C_pC_S =0\,,
\]
which implies $\alpha_0^{\pm}= \pm \sqrt{2C_pC_S} - 2 \sqrt{C_z}$. $\alpha_1$ is determined via the next order equation.
If we consider the traveling wave limit, namely $R_1 \rightarrow \infty$, then clearly $R_2 \rightarrow \alpha_0$ with an {\it algebraic convergence rate}. And if $\alpha_0^+ >0$, or equivalently $C_pC_S > 2C_z$, then $\Omega_2$ persists with width $\alpha_0^+$.

The pressure at $\Gamma_2$ in the limit $R_1 \rightarrow \infty$ simplifies to $\Sigma (R)= \frac{\sqrt{C_z} }{C_S }R_2 +2 \frac{C_z}{C_S}$. When $\alpha_0^+>0$, substituting $R_2$ with $\alpha_0^+$ leads to the following pressure jump:
\[
\Sigma (R)= \sqrt{\frac{2C_zC_p}{C_S}}.
\]

\begin{figure}[h!] 
\centering
\includegraphics[width=0.32\textwidth]{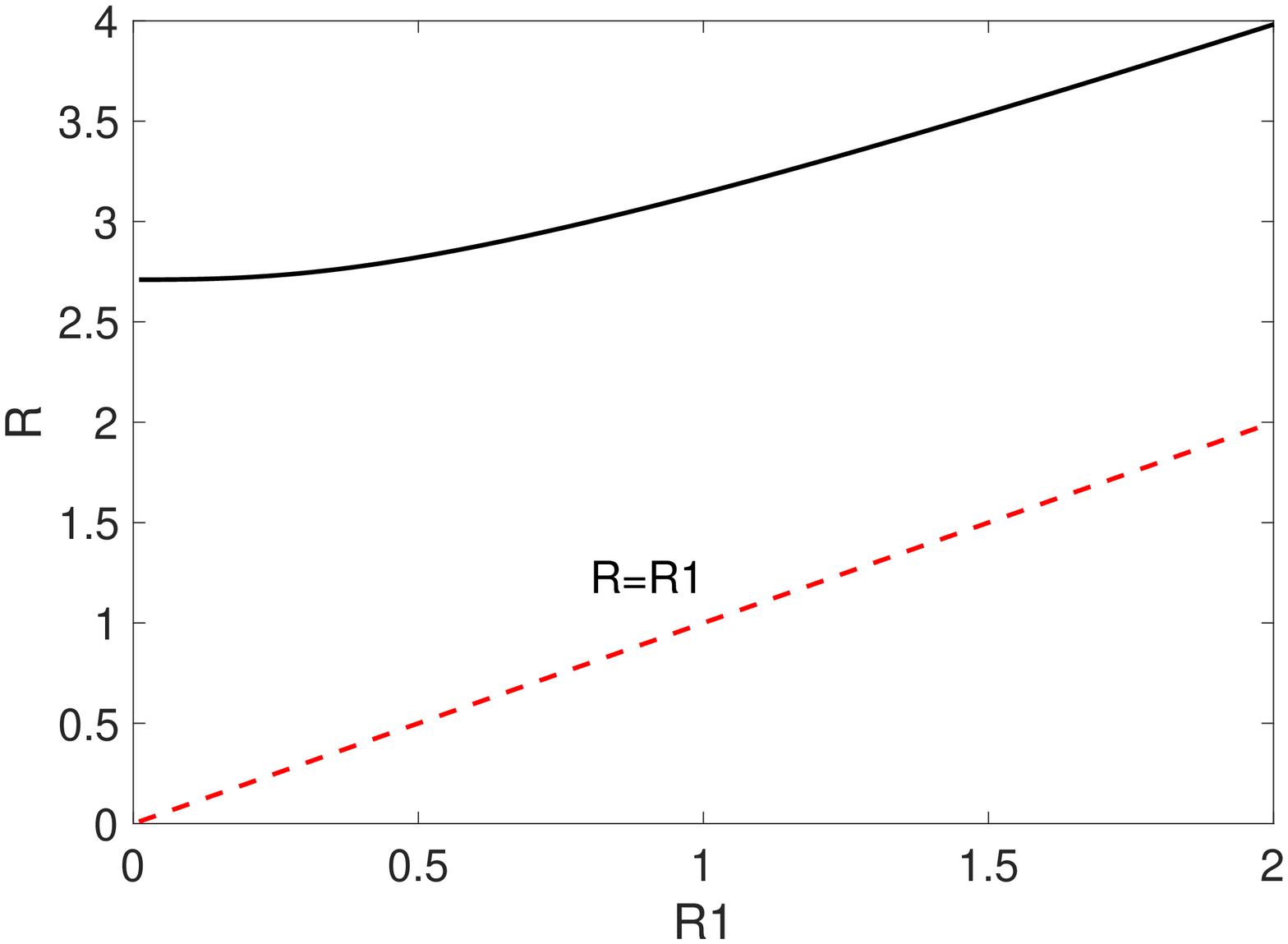}
\includegraphics[width=0.32\textwidth]{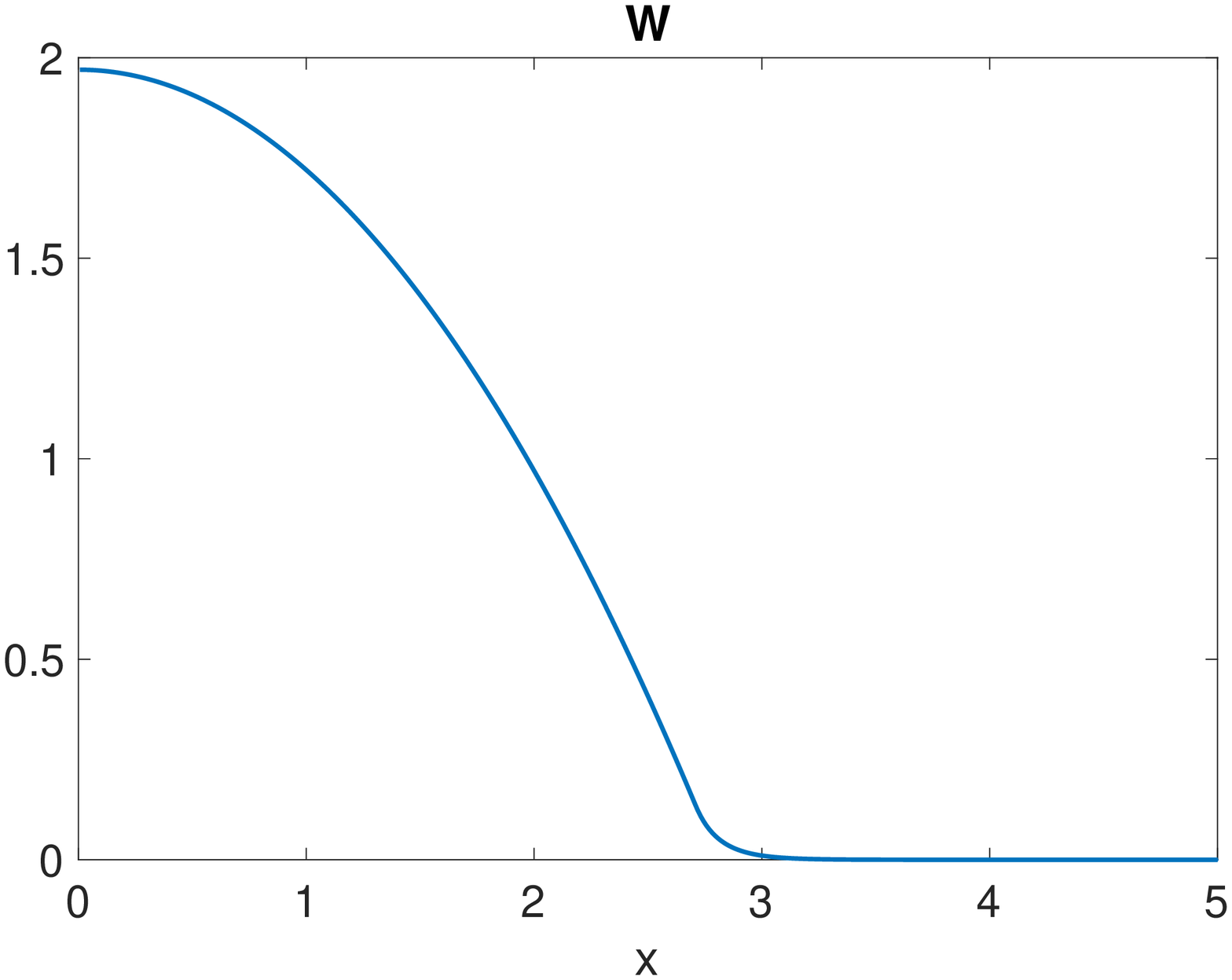}
\includegraphics[width=0.32\textwidth]{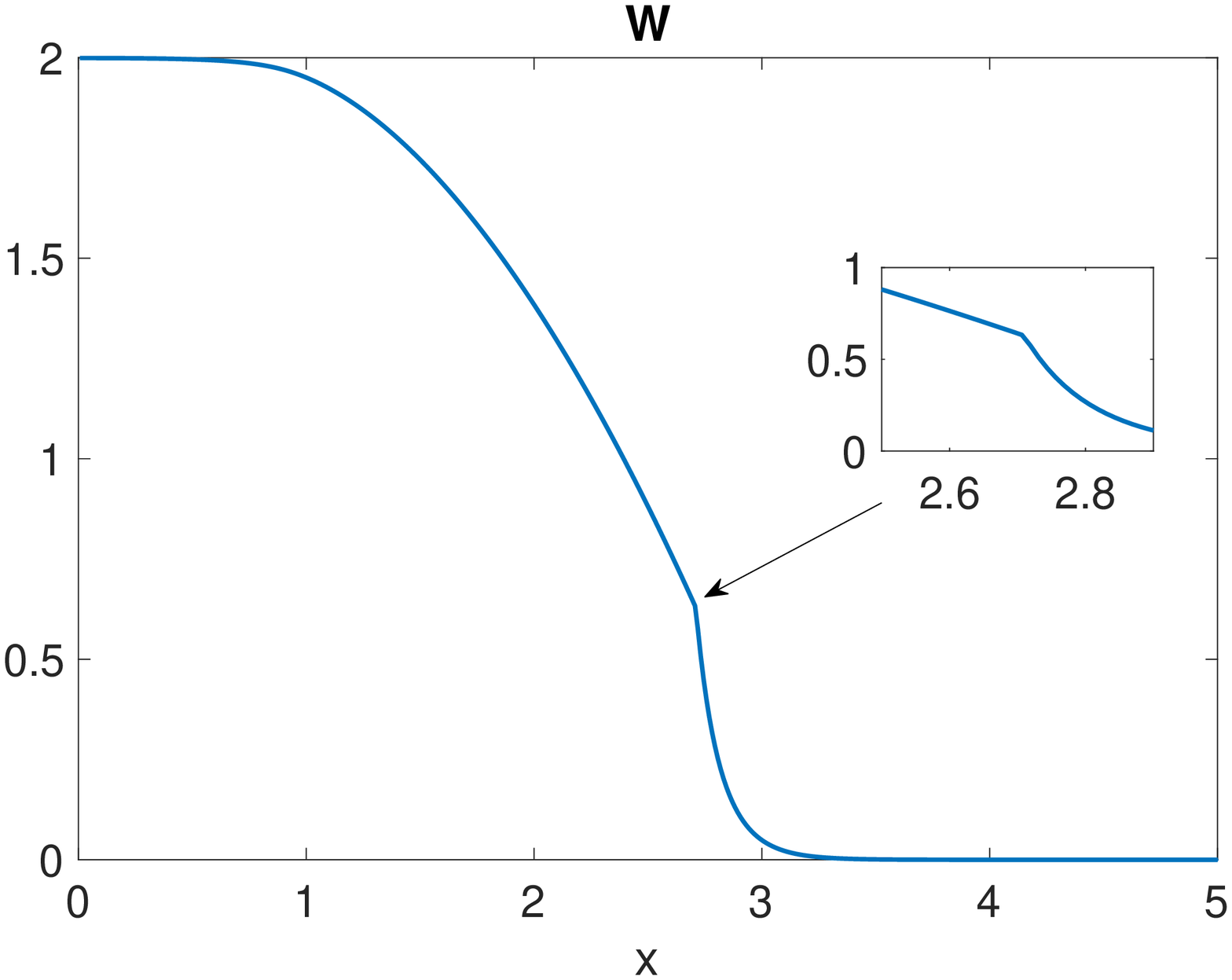}
\caption{Left: relationship between $R_1$ and $R$. Middle: plot of $W(x)$ using \eqref{W-R-2d-lim} with $R = 2.71$ and $R_1 = 0.0151$, satisfying the relation \eqref{rel:2d3}. Right: plot of $W(x)$ using \eqref{W-R-2d-lim} with $R = 2.71$ and $R_1 = 0.9$, not satisfying the relation \eqref{rel:2d3}. 
Here $C_S = 1$,  $C_z = 0.02$,  $C_p = 2$,  $C_{\nu}= 100$, $\eta = 0.01$.}
\label{R-R1-radial}
\end{figure}

When, $R_1$ is finite, there is no explicit solution for $R_2$, instead we consider the evolution equation for $R(t)$
\begin{equation} \label{frel1}
\dot R = - C_S W_r(R) = \frac{R}{2}- \frac{1}{R} \left(  \frac{1}{2}(R_1)^2 - R_1^{} {\sqrt{C_z}}\frac{ I_1 \left(\frac{R_1^{}}{\sqrt{C_z}}\right)}{ I_0 \left(\frac{R_1^{}}{\sqrt{C_z}}\right)} \right).
\end{equation}
Therefore, \eqref{frel1} and \eqref{rel:2d3} can be viewed as a differential-algebraic system of equations, and we can numerically solve for $R_1$ and $R$. Like before, the solvability calls for $R\geq R_1 \geq 0$. We display the relationship between $R$ and $R_1$ from \eqref{rel:2d3} in Figure ~\ref{R-R1-radial}, where a monotone relation is observed.

Further, when $\alpha_0^+>0$, in the traveling wave limit $R_1 \rightarrow \infty$, we obtain
\[
\dot R =\alpha_0^+ + \sqrt{C_z} = \sqrt{2C_pC_S }  - \sqrt{C_z}\,,
\]
which is the same speed as we obtained in 1D case. 

We want to point out that, the significant difference between the 1D and 2D cases is in 2D the effect of curvatures becomes manifest. Indeed, in 2D, when $R_1 \gg 1$, as the asymptotic analysis above shows, the free boundary limit converges to the traveling wave solution only in an algebraic rate. In particular, when $R_1 \gg 1$, we have $R_2 \sim \alpha_0 + \kappa \alpha_1 $, where $\kappa = 1/R$ is the curvature of the tumor front, which asymptotically determines the first order correction of the front propagation speed and the pressure jump. On the contrary, in 1D, the free boundary limit converges to the traveling wave limit exponentially.

\subsection{3D spherical symmetric case}
Similar to the 1D case, with the radial symmetric assumption, we can explicitly solve for the ansatz solution to the regularized impressible model. Here we only list the results, and interested readers can refer to Appendix \ref{3d:cal} for details. The solution to the incompressible limit model takes the following form: 
\begin{equation*}
W(r) = \left\{ \begin{array}{cc} 
C_p + A \, i_0 \left(\frac{r}{\sqrt{C_z}}\right), & r\in \Omega_1
\\ - \frac{1}{6C_S} r^2 + a \frac 1 r +b, &  r \in \Omega_2
\\d  k_0\left(\frac{r}{\sqrt{C_z}} \right), & r \in \Omega_3
\end{array} \right.
\qquad 
\Sigma(r) = \left\{ \begin{array}{cc} 
C_p , & r\in \Omega_1
\\ - \frac{1}{6C_S} r^2 + a \frac 1 r +b + \frac{C_z}{C_S}, &  r \in \Omega_2
\\ 0, & r \in \Omega_3\,,
\end{array} \right.
\end{equation*}
where the parameters are listed below
\begin{align*} 
&a= -\frac{1}{3C_S}(R_1)^3 + (R_1)^2 \frac{\sqrt{C_z}}{C_S} \frac{ i_1 \left(\frac{R_1}{\sqrt{C_z}}\right)}{ i_0 \left(\frac{R_1 }{\sqrt{C_z}}\right)}, \quad 
b=C_p - \frac{C_z}{C_S}+ \frac{(R_1)^2}{2C_S}   -R_1  \frac{\sqrt{ C_z}}{C_S} \frac{ i_1 \left(\frac{R_1}{\sqrt{C_z}}\right)}{ i_0 \left(\frac{R_1 }{\sqrt{ C_z}}\right)},
\\  &A= - \frac{C_z}{C_S i_0 \left(\frac{R_1 }{\sqrt{C_z}}\right) }, \qquad 
d k_0\left(\frac{R}{\sqrt{C_z}} \right)=- \frac{1}{6C_S} R^2 + a \frac{1} {R} +b.
\end{align*}

To get relationship between two boundaries $R$ and $R_1$, again by continuity of $W_r$ at $R$, one has %
\begin{equation}\label{rel:3d2}
 - \frac{d}{\sqrt{C_z}}k_1 \left(\frac{R}{\sqrt{C_z}} \right) = - \frac{1}{3C_S} R - \frac{a}{R^2}. 
\end{equation}
Using $R=R_1+R_2$, it becomes
\begin{equation}\label{rel:3d3}
\sqrt{C_z} \frac{k_0\left(\frac{R}{\sqrt{C_z}} \right)}{k_1 \left(\frac{R}{\sqrt{C_z}} \right)} \left( \frac{1}{3C_S} R+ \frac{a}{R^2} \right) =- \frac{1}{6C_S} R^2 + a \frac{1} {R} +b.
\end{equation}

Now we consider the case when
\[
R_1 \gg 1, \quad R_2=O(1).
\]
By asymptotically expanding each side of \eqref{rel:3d3}, we get
\begin{align*}
\text{L.H.S.} = & \sqrt{C_z} \left(1 - \frac{\sqrt{C_z}}{R_1}  \right)  \left( \frac{1}{C_S}R_2 + \frac{\sqrt{C_z}}{C_S}  + \frac{1}{R_1} \left( -\frac{R_2^2}{C_S}- R_2 \frac{2\sqrt{C_z}}{C_S}  - \frac{C_z}{C_S} \right) \right) + o \left( \frac{1}{R_1} \right).
\end{align*}
Here, we have used the fact that, when $z \gg 1$
\[
\frac{i_1(z)}{i_0(z)}= 1 - \frac{1}{z}+ o\left( \frac 1 z \right), \quad \frac{k_0(z)}{k_1(z)}= 1 - \frac{1}{z}+ o\left( \frac 1 z \right).
\]
Similarly, on the right hand side, we have
\begin{align*}
\text{R.H.S.} 
= &  - \frac{1}{2C_S}R_2^2 - R_2 \frac{\sqrt{C_z}}{C_S}  + C_p - \frac{C_z}{C_S} + \frac{1}{R_1} \left( \frac{R_2^3}{3C_S} + {R_2^2}\frac{\sqrt{C_z}}{C_S}+ \frac{R_2 C_z}{C_S} \right)+ o \left( \frac 1 {R_1} \right).
\end{align*}

To match the terms order by order, we assume when $R_1 \gg 1$,
\[
R_2 = \alpha_0 + \frac{\alpha_1}{R_1} + o\left( \frac 1 {R_1} \right)\,,
\]
then the leading order terms read
\[
(\alpha_0)^2+ 4 \sqrt{C_z}  \alpha_0 +4C_z-2C_pC_S =0.
\]
which implies $\alpha_0^{\pm}= \pm \sqrt{2C_pC_S} - 2 \sqrt{C_z}$.


In the traveling wave limit, namely $R_1 \rightarrow \infty$, then clearly $R_2 \rightarrow \alpha_0$ with an {\it algebraic convergence rate}. If further $\alpha_0^+ >0$, or equivalently $C_pC_S >2C_z$, then $\Omega_2$ persists with width $\alpha_0^+$ in the limit.

As $R_1 \rightarrow \infty$, the pressure at $\Gamma_2$  simplifies to $\Sigma (R)= \frac{\sqrt{C_z} }{C_S }R_2 +2 \frac{C_z}{C_S}$. When $\alpha_0^+>0$, we substitute $R_2$ with $\alpha_0^+$, we obtain the following pressure jump
\[
\Sigma (R)= \sqrt{\frac{2C_zC_p}{C_S}}.
\]

%

Next, we check the front moving speed. 
When, $R_1$ is finite, there is no direct explicit solution for $R_2$. 
Observe that the $R(t)$ satisfies
\begin{equation} \label{frel3}
\dot R = - C_S W_r(R) = \frac{R}{3}+ \frac{1}{R^2} \left(  -\frac{1}{3}(R_1)^3 + (R_1)^{2} {\sqrt{C_z}}\frac{ i_1 \left(\frac{R_1^{}}{\sqrt{C_z}}\right)}{ i_0 \left(\frac{R_1^{}}{\sqrt{C_z}}\right)} \right).
\end{equation}
Thus, \eqref{frel3} and \eqref{rel:3d3} can be viewed as a differential-algebraic system of equations, and we can numerically solve for $R_1$ and $R$ from this system. When  $\alpha_0^+>0$, in the traveling wave limit, $R_1 \rightarrow \infty$, \eqref{frel3} reduces to 
\[
\dot R =\alpha_0^+ + \sqrt{C_z} = \sqrt{2C_pC_S}  - \sqrt{C_z}.
\]

\section{Numerical scheme}

In this section, we introduce a numerical scheme for solving the cell density model \eqref{main}\eqref{eq:Brinkman}. Our goal is to design a scheme that works for a wide range of $C_\nu$ and thus can simulate solutions solutions to the free boundary model when $C_\nu \rightarrow \infty$. 

When $C_\nu$ is large, from the definition of $\Sigma$ in \eqref{constituion0}, the dependence of $\Sigma$ on $\rho$ becomes intractable. More precisely, a small error in $\rho$ induces a big change in $\Sigma$. On the other hand, to find the correct front speed numerically, $\Sigma$ has to be accurate enough. Therefore it is not an easy task of designing numerical schemes that can capture the right solution behavior when $C_\nu$ is large. Other numerical methods developed for degenerate diffusion equation \cite{AWZ,BF,BCW,KRT} only work for $C_\nu$ is $O(1)$. 

Since the incompressible limit \eqref{limit000} is obtained directly from the evolution equation for pressure $\Sigma$ \eqref{eq:Sigma}, we propose a 3-stage prediction-correction-projection method that gives the correct border velocity for $C_{\nu}=O(1)$ and also for $C_{\nu} \gg 1$. This method is essentially inspired by \cite{LTWZ18}, but the prediction-projection object is changed to the potential $W$. In order not to obscure the focus of the current work, we avoid numerical analysis for the method, and save it for future works.

\subsection{The semi-discrete method}

In this part, we introduce the semi-discrete scheme by considering the following augmented system 
\begin{subequations} \label{aug000}
\begin{numcases}{}
\partial_t  \rho - C_S \nabla \cdot (\rho \nabla W) = \Phi (\Sigma, \rho)\,,  \label{aug001}
\\ - C_z \Delta W +W =\Sigma\,, \label{aug002}
\\ \partial_t \Sigma  -C_S \nabla \Sigma \cdot \nabla W - C_S C_{\nu} \Delta W = C_{\nu} H\,,  \label{aug003}
\end{numcases}
\end{subequations}
where $\Sigma$ relates to $\rho$ through the constitution relation \eqref{constituion0}. Note that $\nabla W$ is important in driving $\rho$ forward in time, we combine \eqref{aug002} and \eqref{aug003} to derive the following evolution equation for $W$,
\begin{equation} \label{eq:Wt}
\frac{\partial }{\partial t } \left( W- C_z \Delta W \right) -C_S \nabla \Sigma \cdot \nabla W - C_S C_{\nu} \Delta W = C_{\nu} H.
\end{equation}
Then our semi-discrete predictor-corrector scheme reads as follows. Given $W^n$, $\rho^n$ and $\Sigma^n$, we have 
\begin{subequations} \label{semi000}
\begin{numcases}{}
\frac{(W^* - C_z \Delta W^*)- (W^n-C_z \Delta W^n)}{\Delta t}   -C_S \nabla \Sigma^n \cdot \nabla W^n - C_S C_{\nu} \Delta W^* = C_{\nu} H^n\,, \qquad \label{semi1}
\\ \frac{ \rho^{n+1} - \rho^n}{\Delta t} - C_S \nabla \cdot (\rho^{n+1} \nabla W^*) = \rho^{n+1} H^n\,,   \label{semi2}
\\- C_z \Delta W^{n+1} +W^{n+1} =\Sigma(\rho^{n+1})\,.  \label{semi3}
\end{numcases}
\end{subequations}
where we have used 
\[
\Sigma^n= \Sigma(\rho^n), \quad H^n= H(C_p-\Sigma^n).
\]

Notice that when $C_\nu \rightarrow \infty$,  \eqref{semi1} formally reduces to
\[
- C_S \Delta W^* = H^n,
\]
and therefore captures the free boundary limit. The use of \eqref{semi3} is to dynamically reinforce the constitutive relation between $\rho$ and $W$, which also turns out to be important for stability purpose. 

\subsection{Fully discrete scheme in 1D}
In this part, we elucidate the spatial discretization and form a fully discrete scheme for both one dimensional case.
The one dimensional version of \eqref{semi000} reduces to 
\begin{subequations} \label{semi000-1D}
\begin{numcases}{}
\frac{(W^* - C_z \partial_{xx} W^*)- (W^n-C_z \partial_{xx}  W^n)}{\Delta t}   -C_S \partial_x \Sigma^n   \partial_x W^n - C_S C_{\nu} \partial_{xx}  W^* = C_{\nu} H^n, \qquad \quad  \label{semi1-1D}
\\ \frac{ \rho^{n+1} - \rho^n}{\Delta t} - C_S \partial_{x}   (\rho^{n+1} \partial_{x}  W^*) = \rho^{n+1} H^n\,,   \label{semi2-1D}
\\- C_z \partial_{xx}  W^{n+1} +W^{n+1} =\Sigma(\rho^{n+1})\,.  \label{semi3-1D}
\end{numcases}
\end{subequations}
Then to update $W^*$ from \eqref{semi1-1D}, we have
\[
(1-C_z \partial_{xx} - C_S C_{\nu} \Delta t \partial_{xx}) W^{*} = W^n - C_z \partial_{xx} W^n + \Delta t C_S \partial_x \Sigma^n \partial_x W^n + \Delta t C_{\nu} H^n\,. 
\]
Let $\rho_j$, $\Sigma_j$ and $W_j$ be the approximation of $\rho(x)$, $\Sigma(x)$, and $W(x)$ at position $x_j$, then we approximate the spatial derivatives in the above equation via center difference, i.e., 
\[
(\partial_{xx} W)_j = \frac{W_{j-1} -2 W_j + W_{j+1}}{\Delta x^2} , \quad
(\partial_x W)_j = \frac{W_{j+1} - W_{j-1}}{2\Delta x} , \quad 
(\partial_x \Sigma)_j = \frac{\Sigma_{j+1} - \Sigma_{j-1}}{2\Delta x},
\]
and we use zero boundary condition for both $\Sigma$ and $W$. 

To propagate $\rho$ in time, we use central scheme to treat the convection term in \eqref{semi2-1D}. More specifically, let 
\[
u_{j+\half}^n = C_S\frac{W_{j+1}^* - W_j^* }{\Delta x},
\]
then \eqref{semi2-1D} is discretized as
\begin{align} \label{rho-fullD}
(1+\Delta t H^n)\rho_j^{n+1}  = \rho_j^n + \frac{\Delta t}{2 \Delta x} & \left[  u_{j+\half}^n (\rho_j^{n,R} +  \rho_{j+1}^{n,L}) -|u_{j+\half}^n| (\rho_j^{n,R} -  \rho_{j+1}^{n,L}) \right.
\\ & \quad \left.  - u_{j-\half}^n (\rho_{j-1}^{n,R} +  \rho_{j}^{n,L}) + |u_{j-\half}^n|  (\rho_{j-1}^{n,R} -  \rho_{j}^{n,L}) \right] \,,
\end{align}
where 
\[
\rho_j^{n,R} = \rho_j^n  + \half \sigma_j^n , \quad
\rho_j^{n,L} = \rho_j^n  - \half \sigma_j^n \,,
\]
and
\[
\sigma_j^n = \left\{ \begin{array}{cc} 0 & \text{ if \quad } (\rho_{j+1}^n - \rho_j^n) (\rho_{j}^n - \rho_{j-1}^n) <0 
\\ \rho_j ^n- \rho_{j-1}^n  & \text{else if \quad } |\rho_{j+1}^n - \rho_j^n| > |\rho_j^n - \rho_{j-1}^n|
\\ \rho_{j+1}^n - \rho_{j}^n & \text{else} 
\end{array}. \right.
\]

\subsection{Fully discrete 2D radial symmetric case}
In the same line of \eqref{aug000}, we first write down the augmented system for the two dimensional radial symmetric case
\begin{equation*} 
\begin{cases}{}
\partial_t  \rho - C_S \frac{1}{r} \partial_r (r \rho \partial_r W)= \Phi (\Sigma, \rho)\,,  
\\ - C_z \frac{1}{r} \partial_r (r \partial_rW) +W =\Sigma\,,
\\ \partial_t \Sigma  -C_S \partial_r \Sigma \partial_r W - C_S C_{\nu} \frac{1}{r} \partial_r (r \partial_r W) = C_{\nu} H\,,  
\end{cases}
\end{equation*}
and semi-discrete it in the same manner as in \eqref{semi000-1D} to get 
\begin{subequations} \label{semi000-R}
\begin{numcases}{}
\frac{(W^* - C_z \frac{1}{r} \partial_r ( r \partial_r W^*) )- (W^n-C_z \frac{1}{r} \partial_r ( r \partial_r W^n))}{\Delta t}  \nonumber
\\  \hspace{5cm}-C_S \partial_r \Sigma^n   \partial_r W^n - C_S C_{\nu} \frac{1}{r} \partial_r ( r \partial_r W^*) = C_{\nu} H^n, \qquad \quad  \label{semi1-R}
\\ \frac{ \rho^{n+1} - \rho^n}{\Delta t} - C_S \frac{1}{r} \partial_r (r \rho^{n+1} \partial_r W^*) = \rho^{n+1} H^n\,,   \label{semi2-R}
\\- C_z \frac{1}{r} \partial_r(r \partial_r W^{n+1}) +W^{n+1} =\Sigma(\rho^{n+1})\,.  \label{semi3-R}
\end{numcases}
\end{subequations}
To discretize in space, let $[0,L_r]$ be our computational domain, and denote $r_j = \frac{\Delta r}{2} + (j-1)\Delta r$, $j = 1, 2, \cdots N_r$, then $\rho_j$, $\Sigma_j$ and $W_j$ approximate $\rho(x)$, $\Sigma(x)$, and $W(x)$ at position $r_j$ respectively. Then the spatial discretization for $W$ reads as 
\begin{align*}
\left[\frac{1}{r} \partial_r(r \partial_r W^{n+1}) \right]_j &= \frac{1}{r_j} \left[ (r\partial_r W)_{j+\half} - (r\partial_r W)_{j-\half}\right] \frac{1}{\Delta r} 
\\& = \frac{1}{r_j} \frac{1}{\Delta r^2} \left[  r_{j+\half} W_{j+1}  \!-\! (r_{j+\half} + r_{j-\half}) W_j \!+\! r_{j-\half} W_{j-1} \right]\,~ j = 2, \cdots, N_r-1\,,
\end{align*}
and Neumann boundary condition implies 
\begin{align*}
&\left[\frac{1}{r} \partial_r(r \partial_r W^{n+1}) \right]_1 = 
\frac{1}{r_1} \left[ (r\partial_r W)_{\frac{3}{2}} - (r\partial_r W)_{\half}\right] \frac{1}{\Delta r} 
= \frac{1}{r_1} \frac{r_{\frac{3}{2}} }{\Delta r^2}  \left(  W_2 - W_1 \right)\,,
\\
&\left[\frac{1}{r} \partial_r(r \partial_r W^{n+1}) \right]_{N} = 
\frac{1}{r_N} \left[ (r\partial_r W)_{N+\half} - (r\partial_r W)_{N-\half}\right] \frac{1}{\Delta r} 
= \frac{1}{r_N} \frac{r_{N-\half} }{\Delta r^2}  \left(  W_{N-1} - W_N \right)\,.
\end{align*}
Likewise 
\begin{align*}
(\partial_r W)_j = \frac{1}{2\Delta r} (W_{j+1} - W_{j-1}), \quad j = 2, \cdots, N-1\,,
\end{align*}
and 
\[
(\partial_r W)_1 = \frac{1}{2\Delta r} (W_2 - W_1), \quad 
(\partial_r W)_N = \frac{1}{2\Delta r} (W_N - W_{N-1})\,.
\]
To update $\rho$, let $g(t,r) = r \rho(t,r)$ and $u = C_S \partial_r W$, then \eqref{semi2-1D} is reformulated as
\[
\frac{g^{n+1} - g^n}{\Delta t} - \partial_r ( g^{n+1} u^*)  = g^{n+1} H^n\,.
\]
Denote $u_{j+\half}^n = C_S\frac{W_{j+1}^* - W_j^* }{\Delta r}$, the above equation can be discretized via central scheme similarly as in \eqref{rho-fullD}. 


\section{Numerical results}
In this section, we conduct a few numerical tests to verify both the model consistency and numerical scheme's efficiency and accuracy. More specifically, we have two major objectives. First, we numerically verify that the cell density model effectively captures the incompressible limit when $C_{\nu} \gg 1$, and in particular, the consistences in the boundary propagation speed and in the pressure jump are carefully checked. Secondly, we compare the boundary moving speeds of the solutions with finite radius with the those of the traveling wave solutions in 1D and 2D respectively. We shall see the difference in the convergence trends between the 1D tests and the 2D tests, which confirms our asymptotic analysis results in Section~\ref{sec:fbm}.


\subsection{1D case}
\begin{figure}[h!] 
\includegraphics[width=0.32\textwidth]{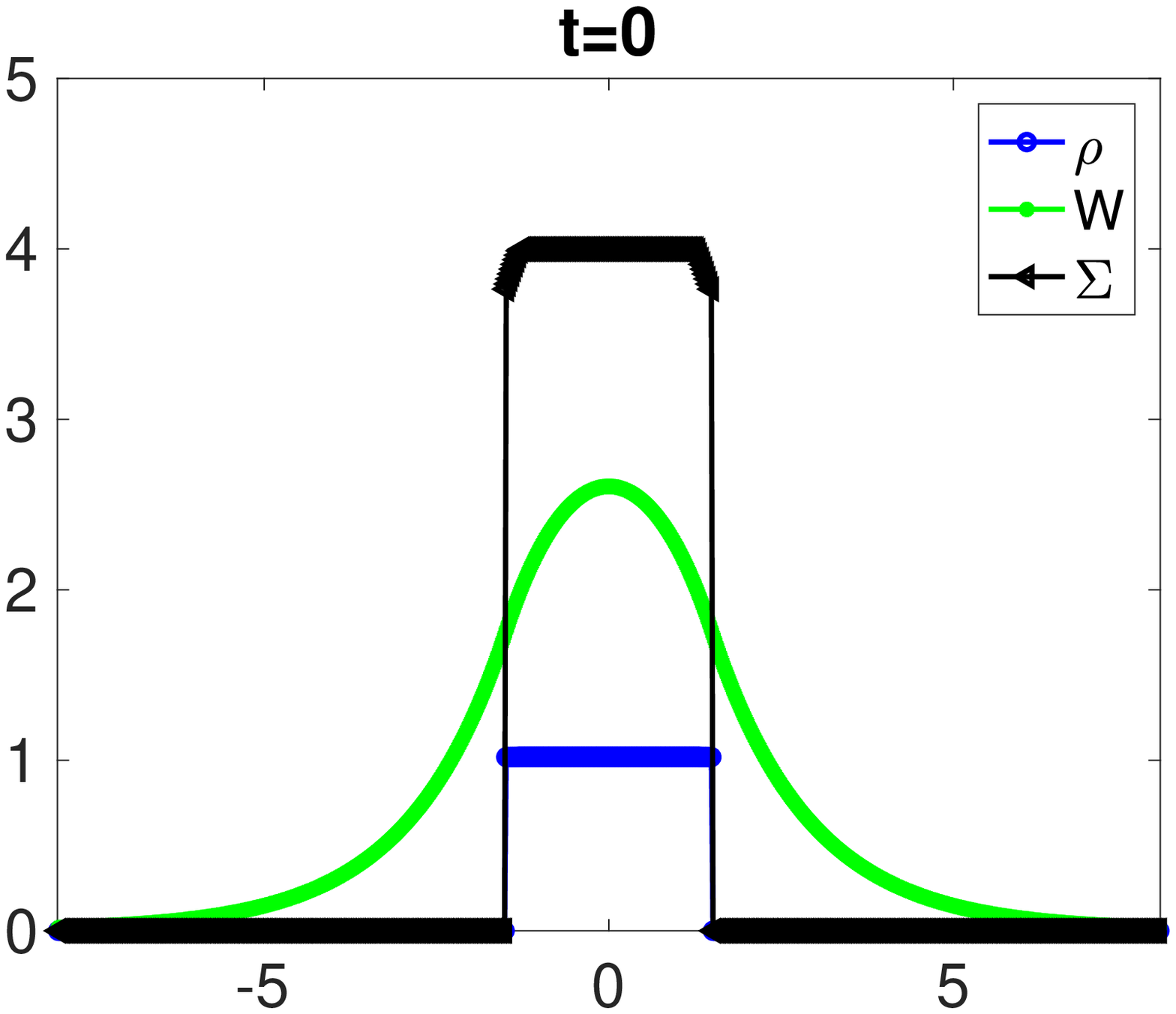}
\includegraphics[width=0.32\textwidth]{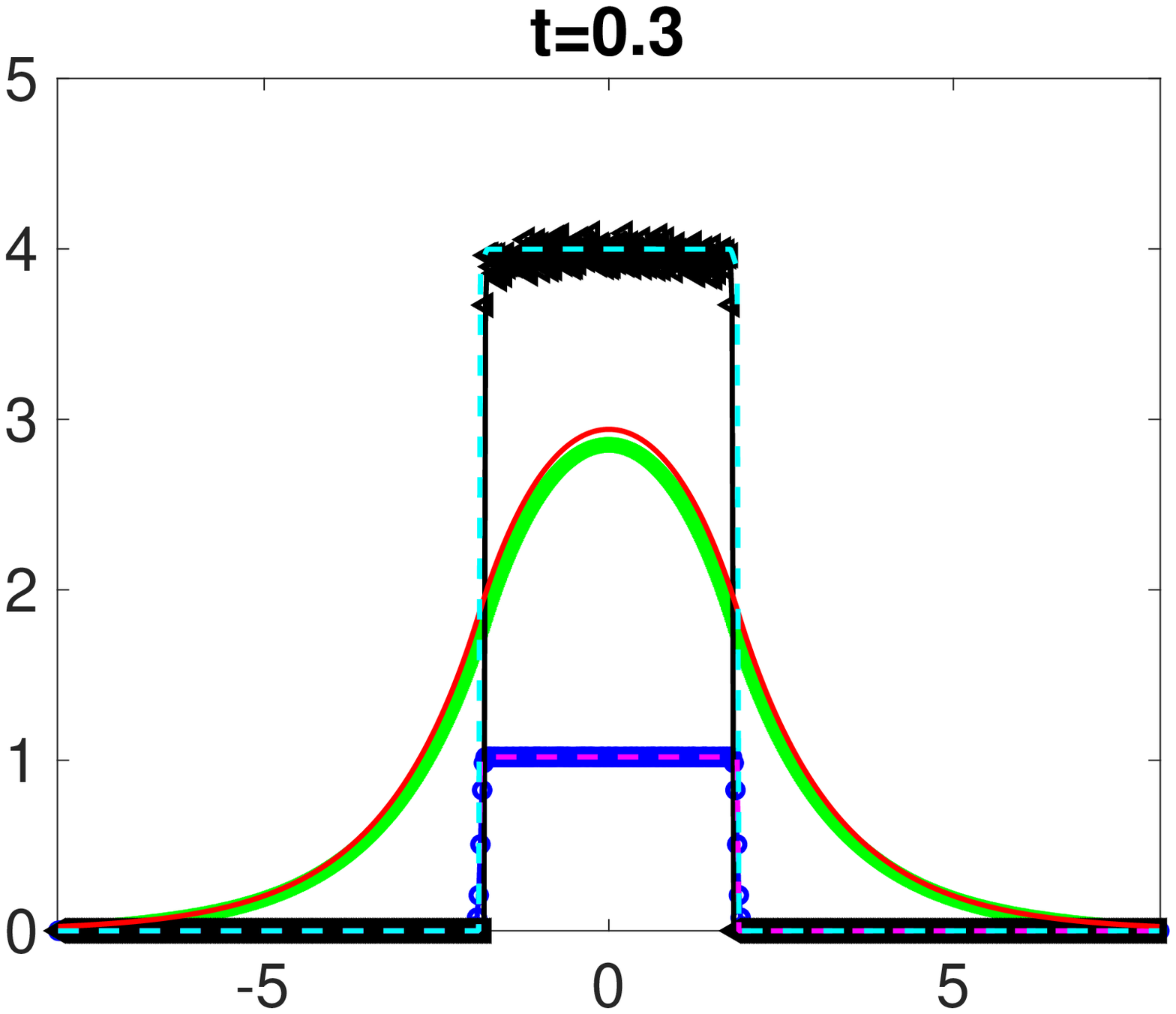}
\includegraphics[width=0.32\textwidth]{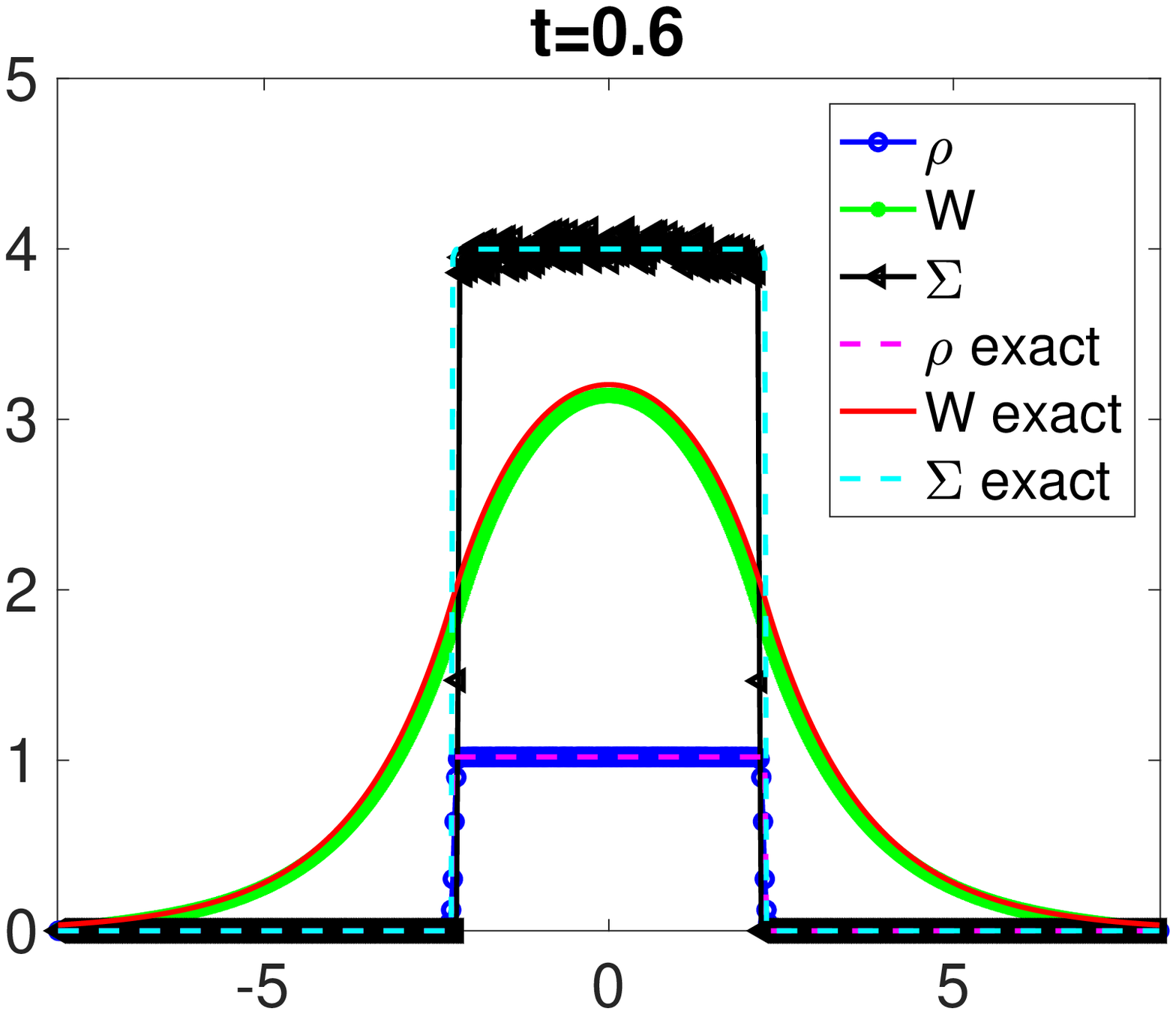}
\caption{1D evolution of $\rho$, $\Sigma$ and $W$ with initial data given by the limit model. Here $\Delta x = 0.025$, $\Delta t = 0.0011$. The parameters used in this example are $C_z = 2$, $C_S = 1$, $C_p = 4$, $C_{\nu} = 200$ and $\eta = 0.001$.}
\label{fig:ex1}
\end{figure}

First, we check the asymptotic property of the scheme \eqref{semi000-1D} when $C_{\nu}$ is sufficiently large. Let $R(0) = 1.5$, choose $R_2(0)$ such that it satisfies \eqref{eq:R2} with $R_1(0) = R(0) - R_2(0)$. Then initial condition $\Sigma(0,x)$ and $W(0,x)$ are chosen of the form \eqref{Sigma-1D} and \eqref{W-1D} respectively, where $\Omega_1(t) = [- R_1(0), R_1(0)]$ and $\Omega_2 = [-R(0), -R_1(0)] \cup [R_1(0), R(0)]$. The constants are $C_z = 2$, $C_S = 1$, $C_p = 4$, $C_{\nu} = 200$. The regularization parameter $\eta = 0.001$. We plot the solution in Figure \ref{fig:ex1}, where a good match between numerical solution to model \eqref{main} \eqref{eq:Brinkman} and exact solution to the limit model \eqref{limit000} is observed. Here the oscillation in $\Sigma$ is due to the amplification by $log$ of the small oscillation in $\rho$ near the interface. 

\begin{figure}[h!] 
\includegraphics[width=0.32\textwidth]{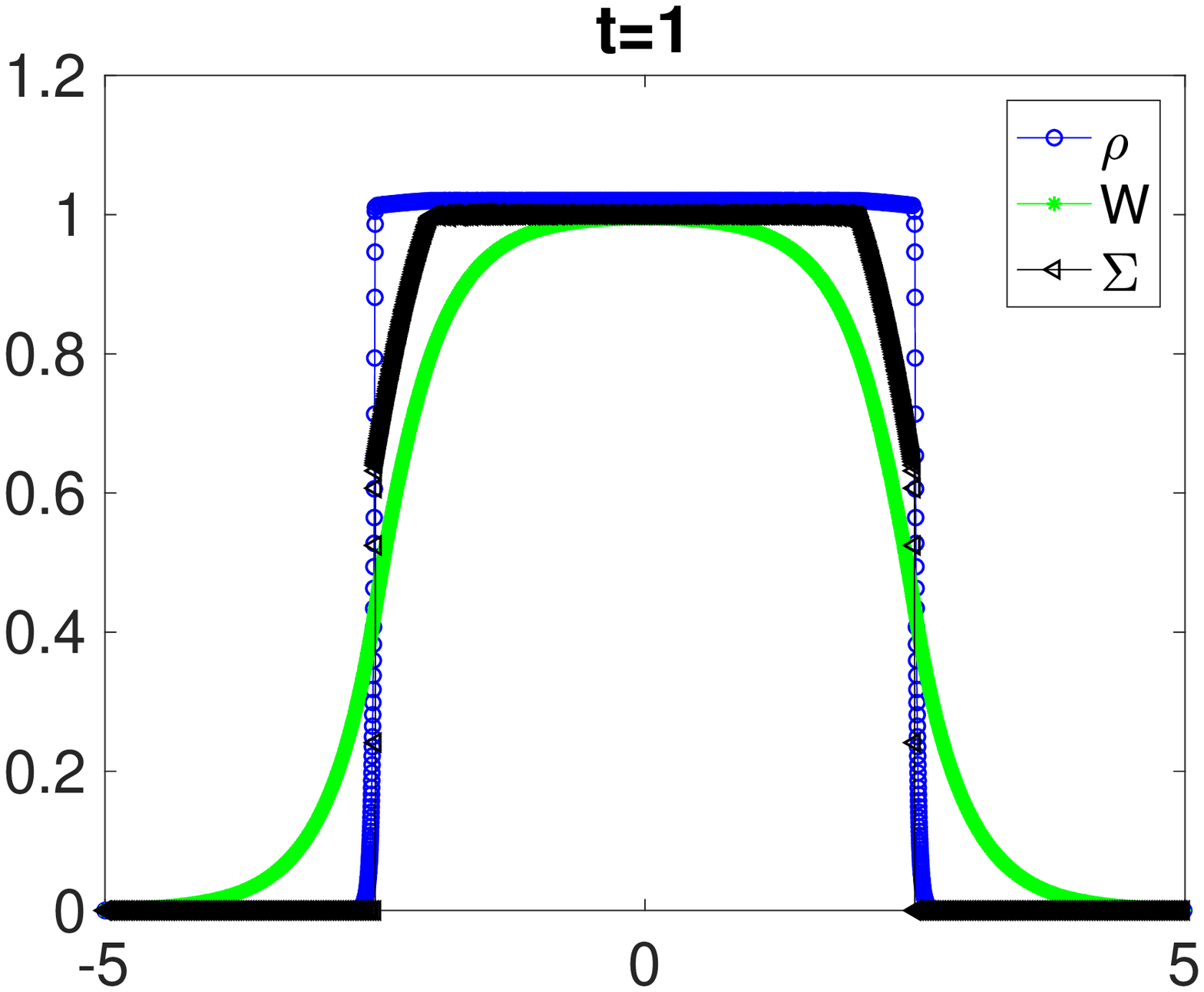}
\includegraphics[width=0.32\textwidth]{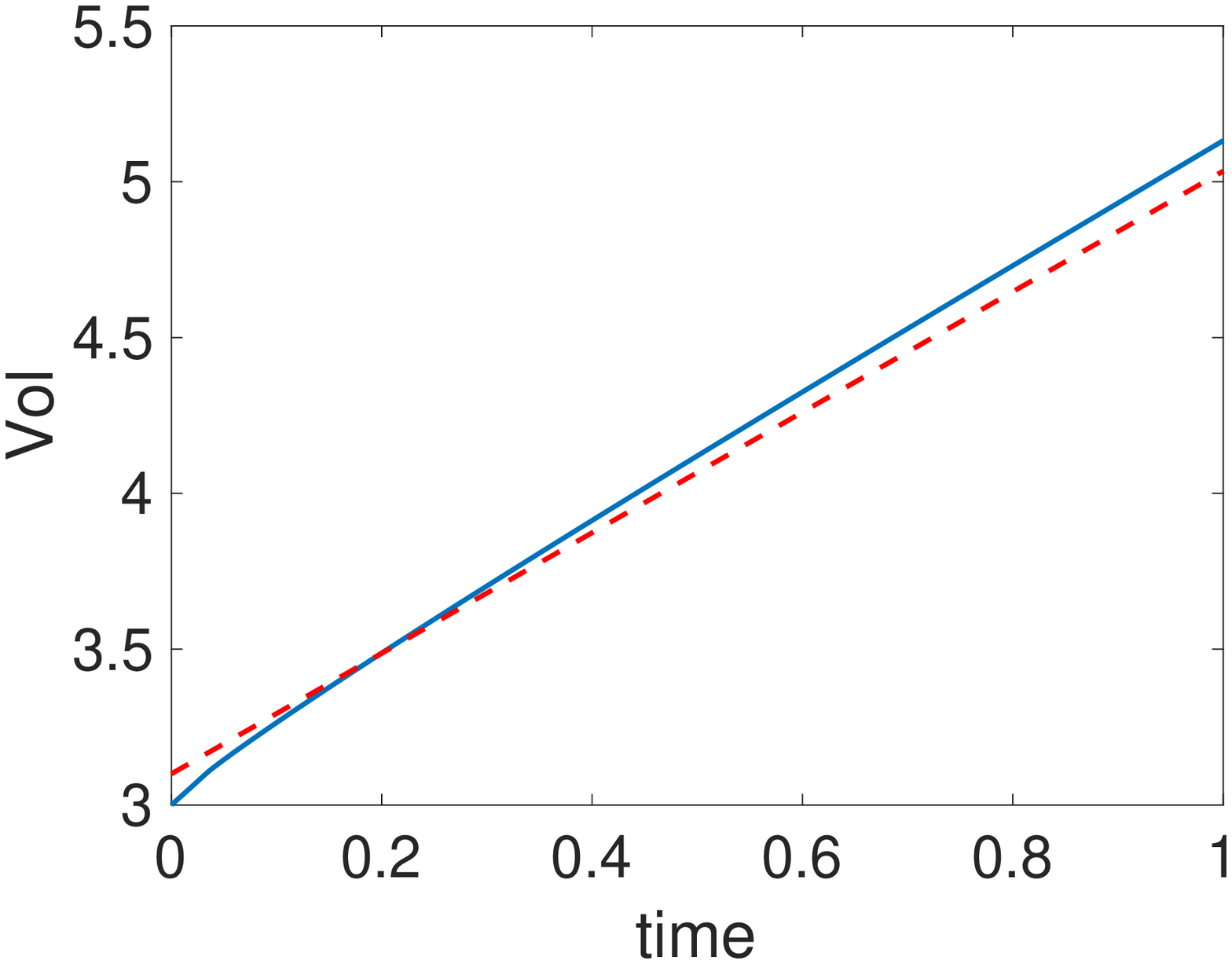}
\includegraphics[width=0.32\textwidth]{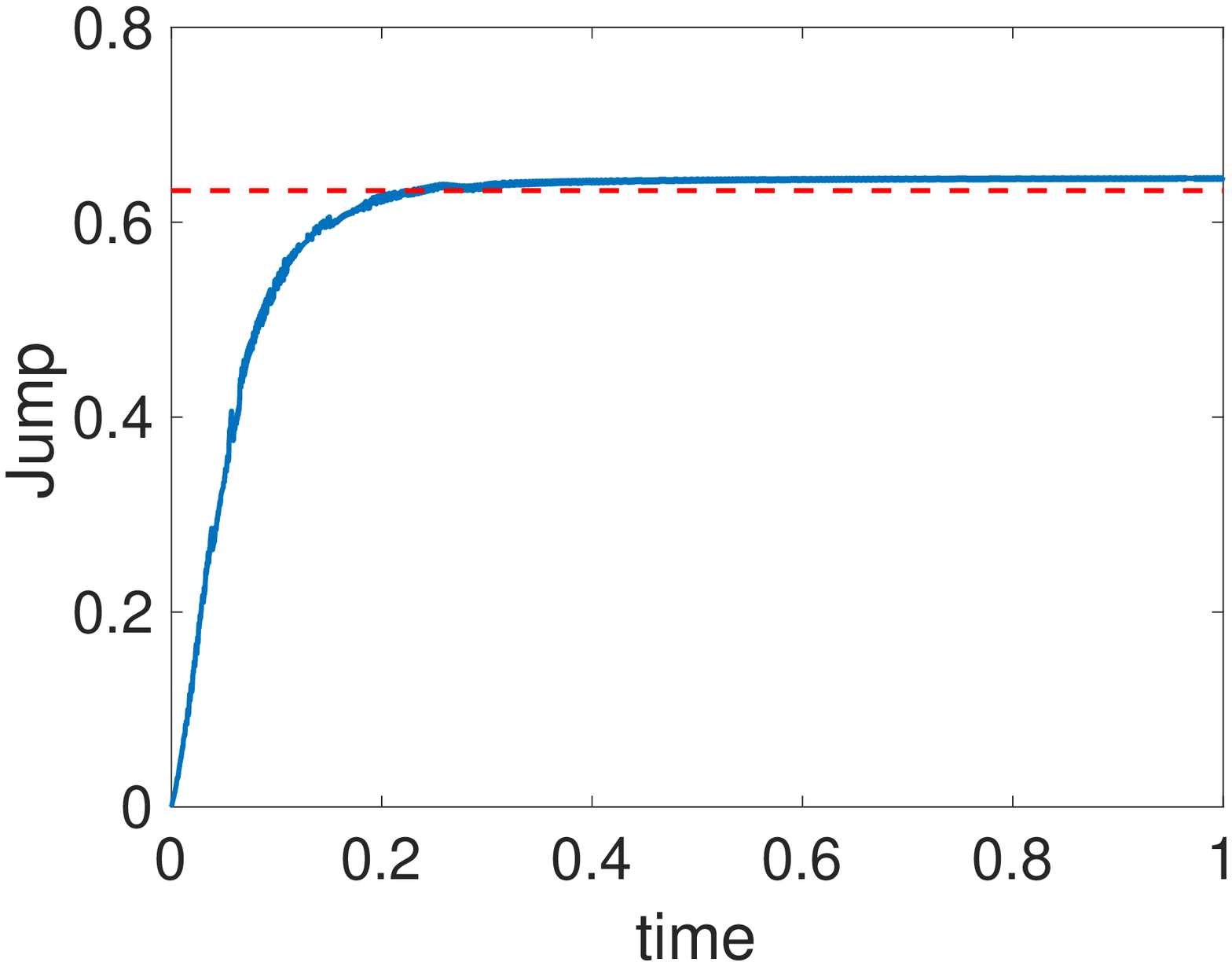}
\caption{1D evolution of cell density model with initial data given by the limit model. Here $\Delta x = 0.0125$, $\Delta t = 9.9829e-05$. Left: plot of $\rho$, $\Sigma$ and $W$ at time $t=1$.  Middle: time versus volume. Here volume is $\int \rho \rd x$ and the red dashed line represents slope $2(\sqrt{2C_pC_s} - \sqrt{C_z})$ given by the analytical formula. Right: jump versus time. The red dotted line is the limit pressure jump $\sqrt{2C_z C_p/C_S}$.}
\label{fig:ex2}
\end{figure}

Next, we check the behavior of the jump in $\Sigma$, the volume of the tumor, and the tumor invading front, versus time. The initial data is again chosen to be of the form \eqref{Sigma-1D} and \eqref{W-1D}  but with $R_1 = 1$, $R = 1.5$. The parameters are $C_z = 0.2$, $C_S = 1$, $C_p = 1$, $C_{\nu} = 50$, and $\eta = 0.001$, and the results are gathered in Figure \ref{fig:ex2}, again good agreement between numerical solution and theoretical predictions are observed. 

\begin{figure}[h!]
\includegraphics[width=0.45\textwidth]{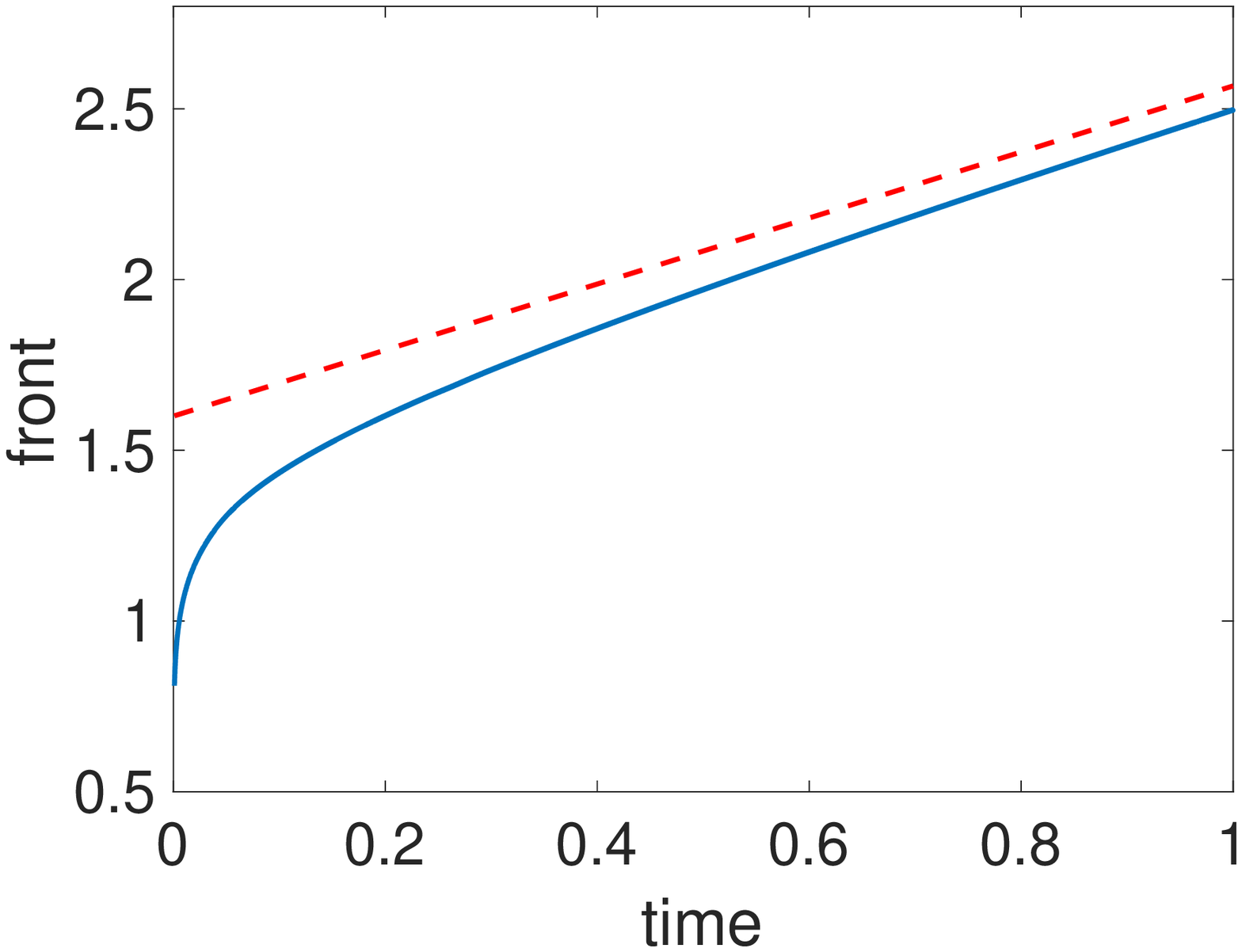}
\includegraphics[width=0.45\textwidth]{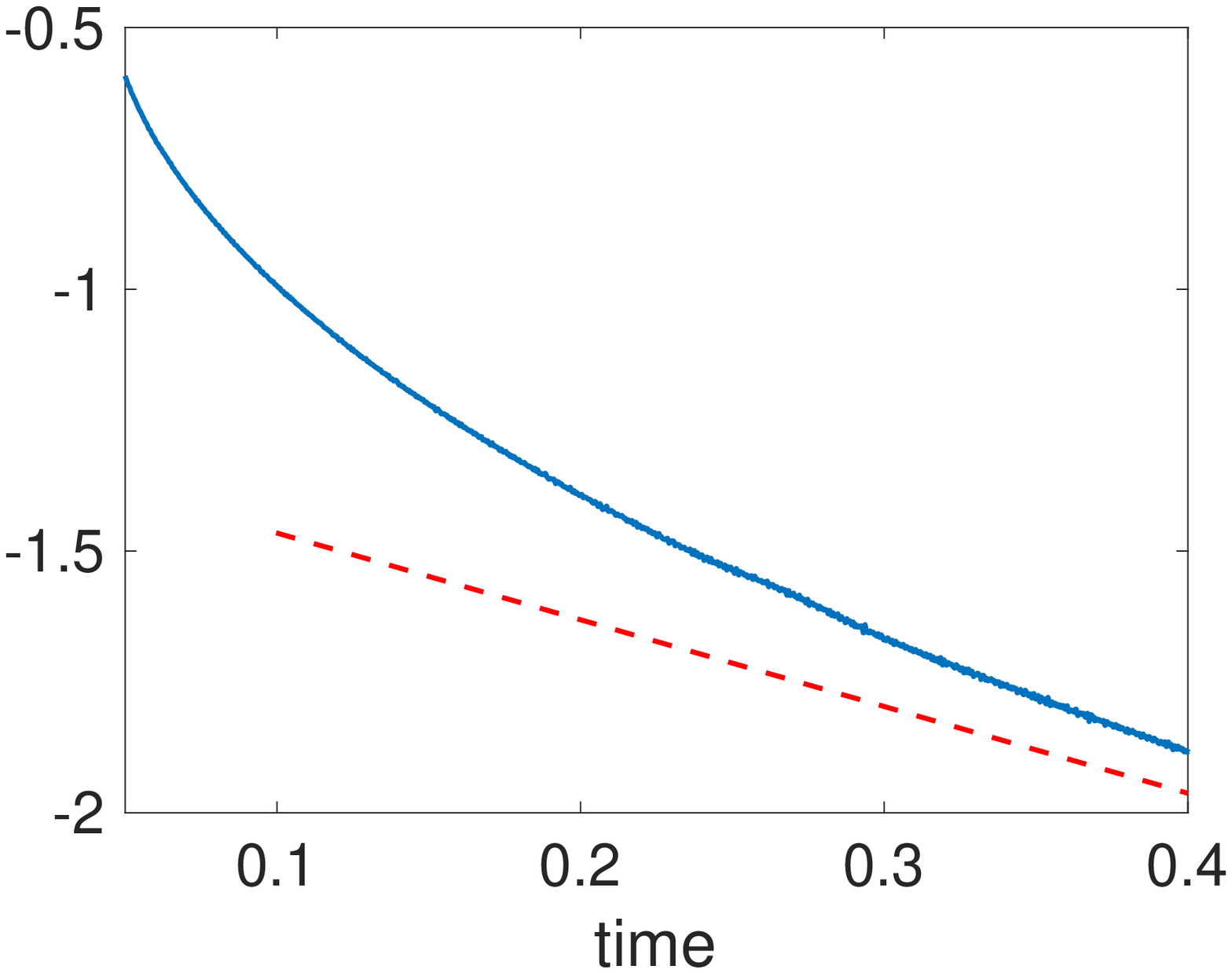}
\caption{Left: plot of the front versus time. The red dashed line indicates slope $\sqrt{2C_p C_S} - \sqrt{C_z}$. Right: $\log|\dot R(t) - (\sqrt{2C_p C_S} - \sqrt{C_z})|$ versus time, where $\dot R(t)$ denotes the front speed; the red dashed line denotes the slope $4(\sqrt{2C_p C_S} - \sqrt{C_z})/\sqrt{C_z}$, which is indicated by the analytical formula.  Here $\Delta x = 0.0125$, $\Delta t = 9.9829e-05$.}
\label{fig:1D22}
\end{figure}

We further check the convergence of propagation speed towards the limit in Figure \ref{fig:1D22}. Here on the left, the dashed line is with slope denoted by the constant speed in the large $R_1$ limit, corresponding to the traveling wave models. One sees that the blue curve, obtained by evolving the cell density model, approaches the red dashed line, indicating that it is the correct asymptote. On the right, an exponential convergence towards the asymptote is displayed.

\subsection{2D radial symmetric case}
\begin{figure}[h!] 
\includegraphics[width=0.32\textwidth]{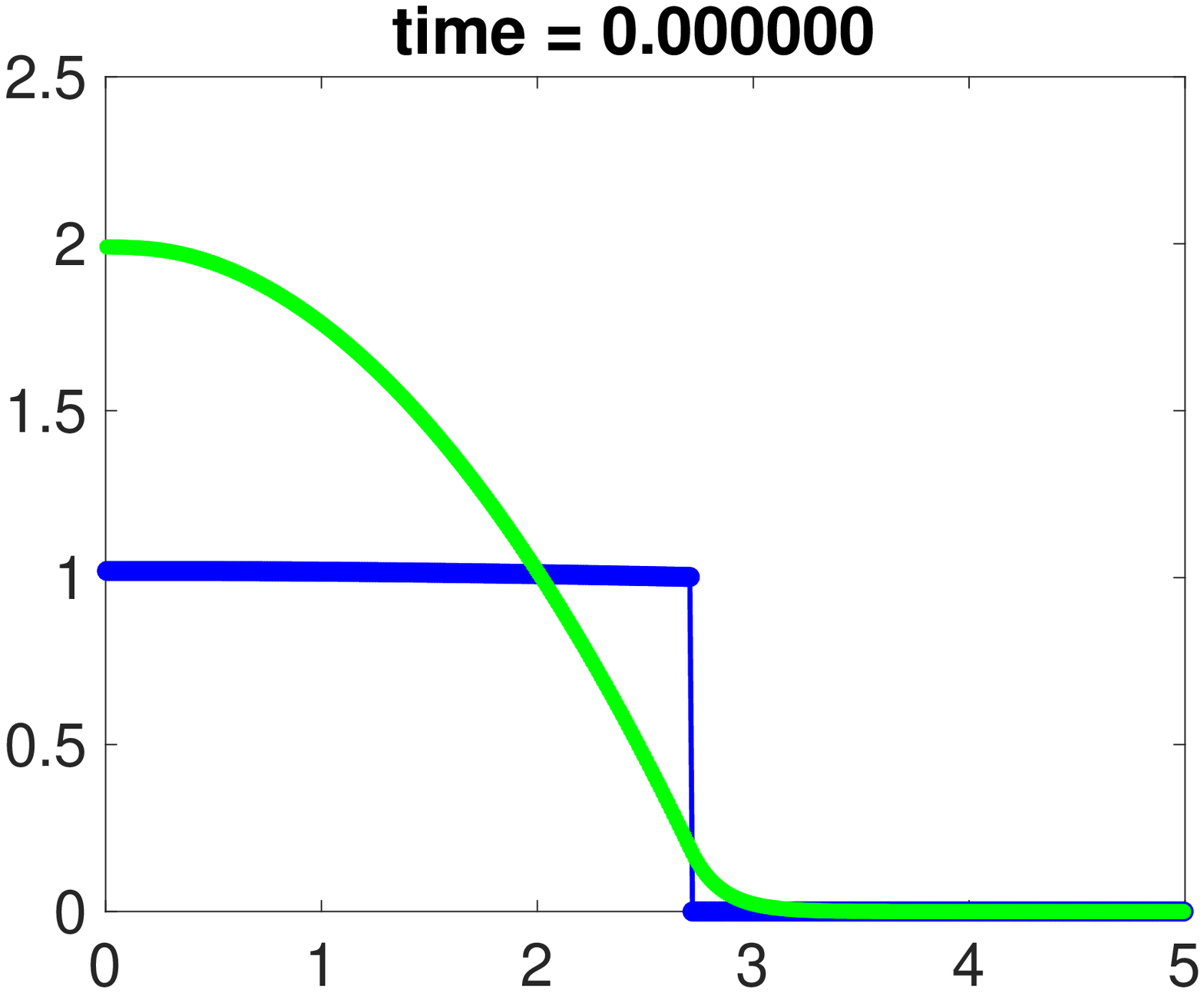}
\includegraphics[width=0.32\textwidth]{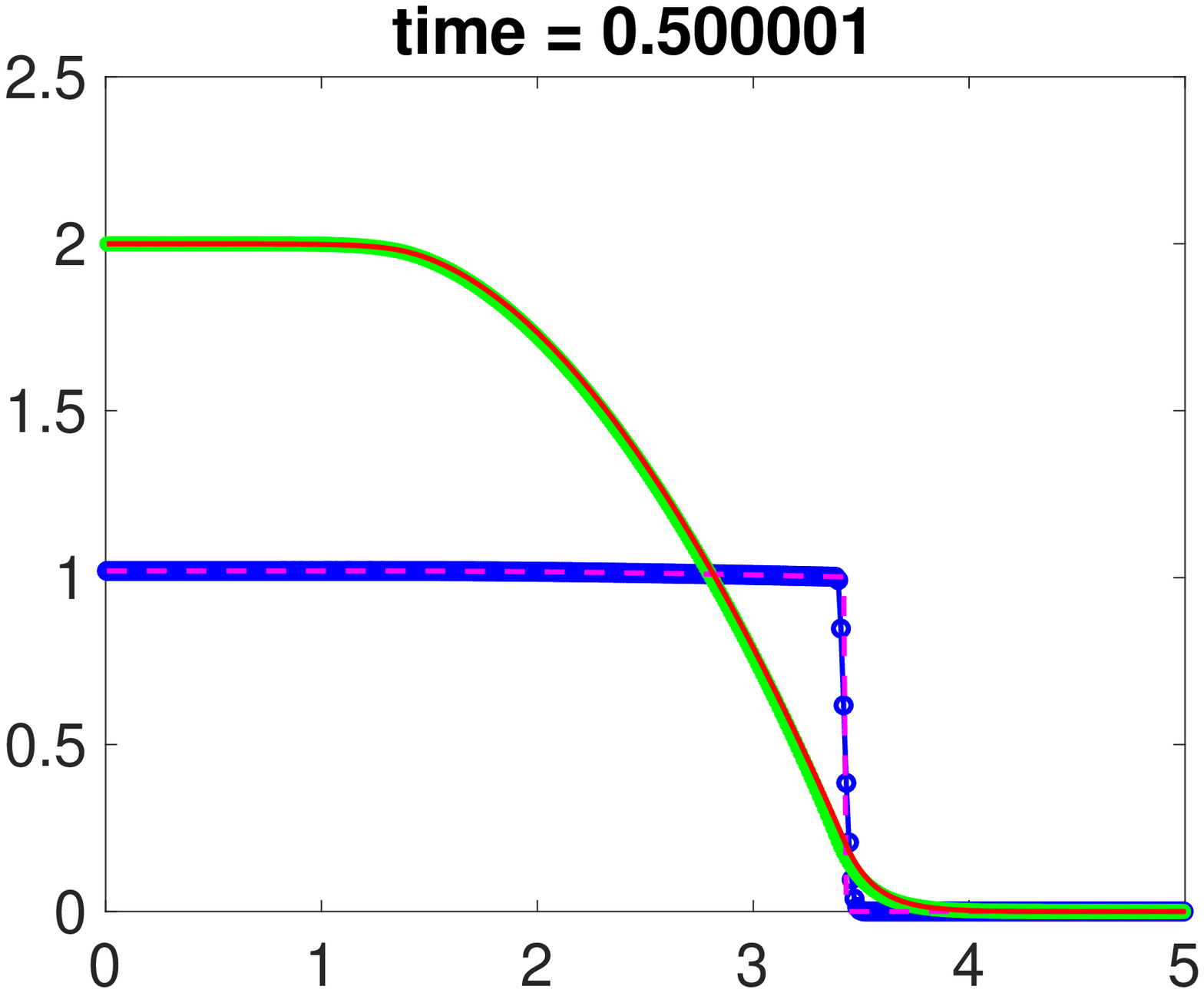}
\includegraphics[width=0.32\textwidth]{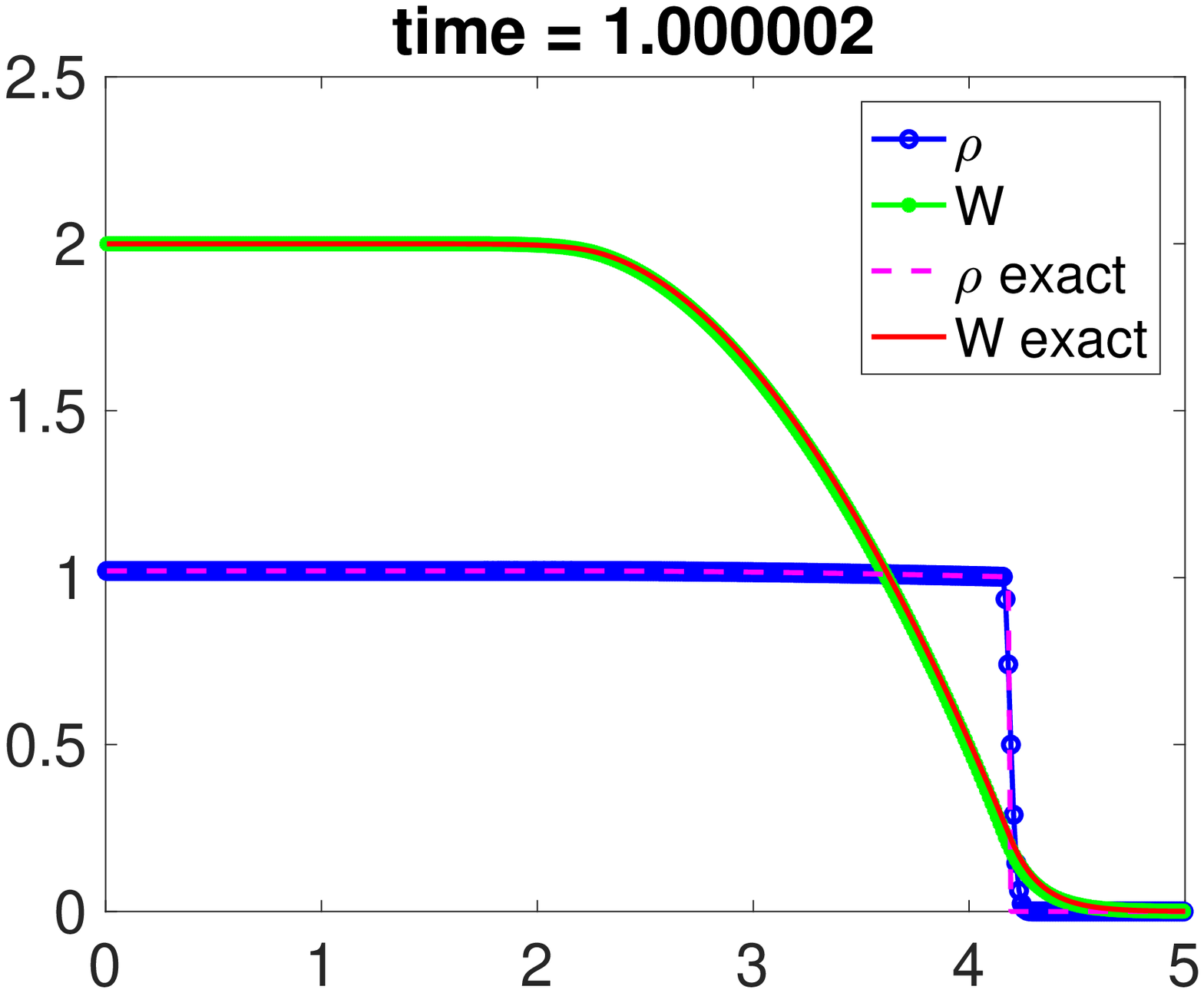}\\
\includegraphics[width=0.32\textwidth]{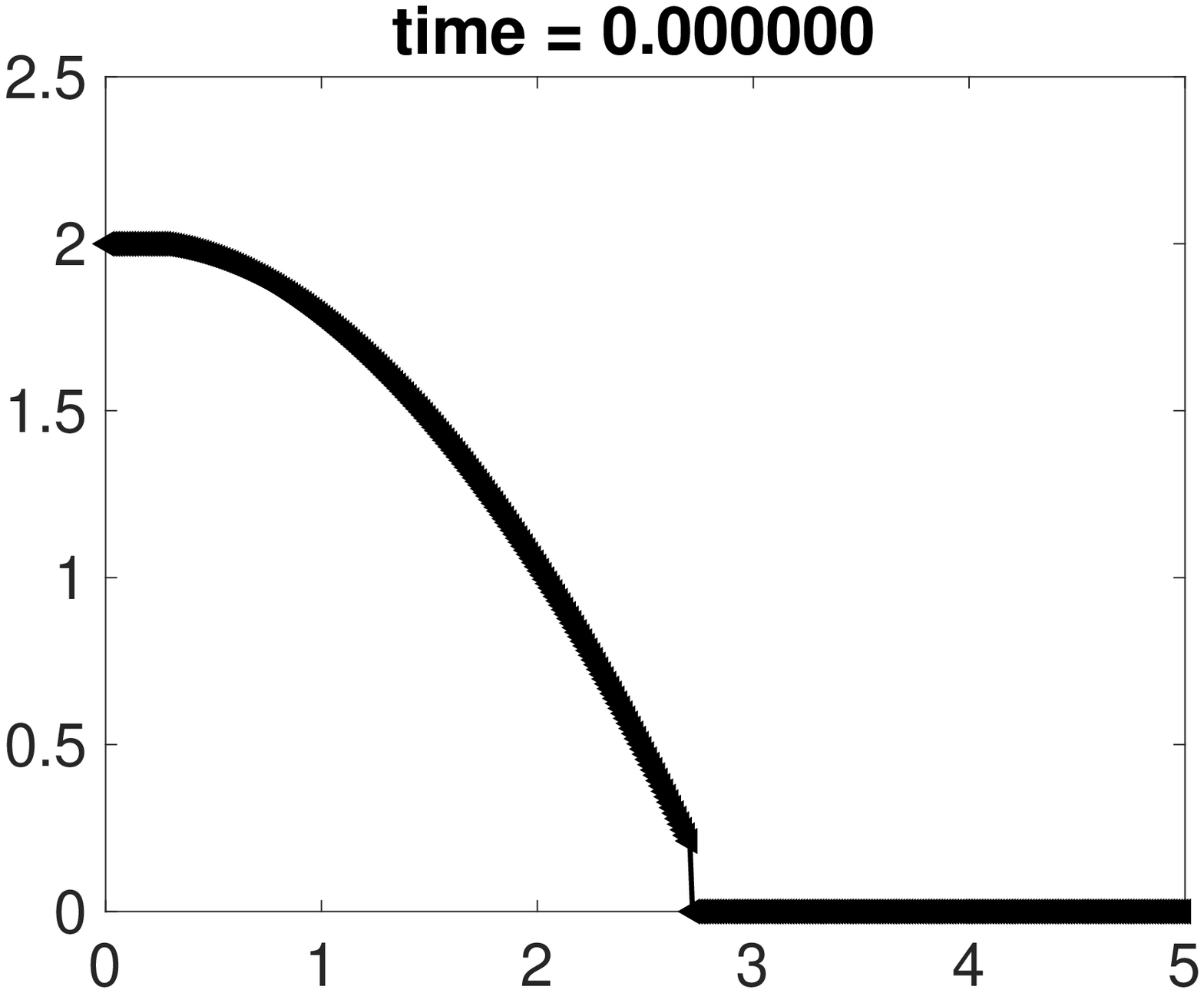}
\includegraphics[width=0.32\textwidth]{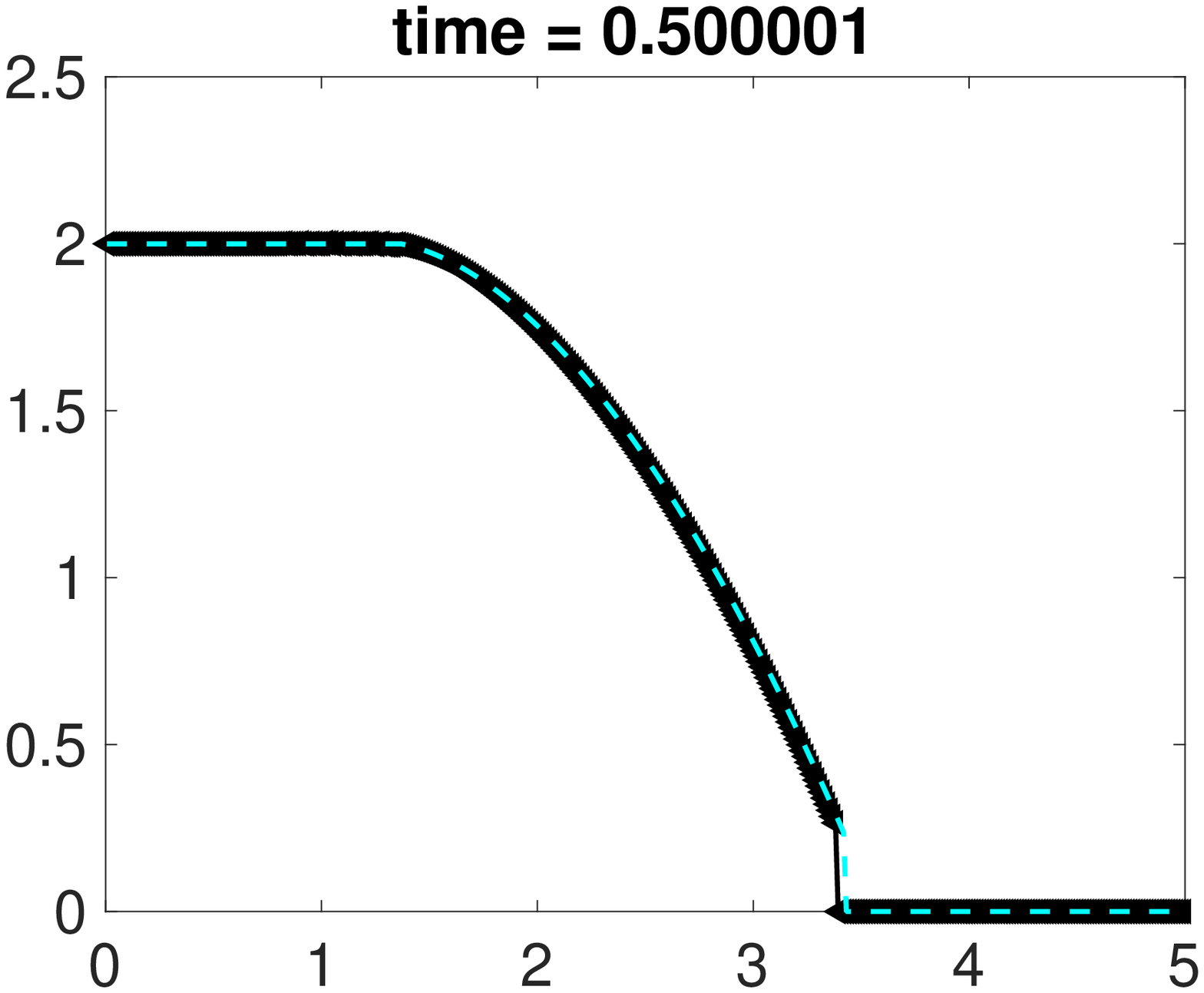}
\includegraphics[width=0.32\textwidth]{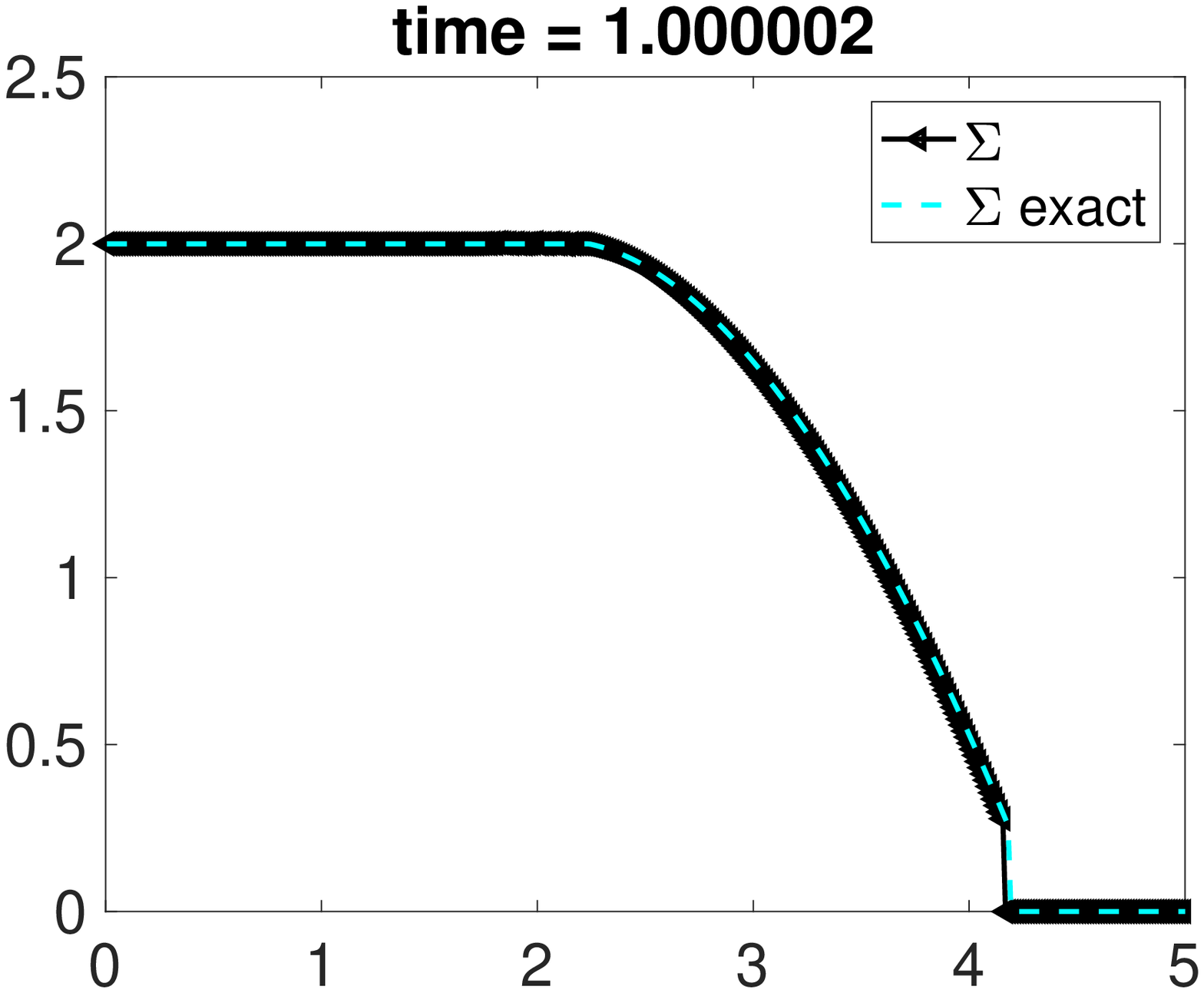}
\caption{2D radially symmetric case. Evolution of $\rho$, $\Sigma$ and $W$ with initial data given by the limit model. Here $\Delta r = 0.0125$, $\Delta t = 5.3351e-5$. The parameters used in this example are $C_{\nu} = 100$, $C_p = 2$, $C_z = 0.02$, $C_s = 1$ and $\eta = 0.0001$.}
\label{fig:radial1}
\end{figure}

In this subsection, we consider 2D radial symmetric case. Like before, we first check the asymptotic property of the scheme \eqref{semi000-R} with sufficiently large $C_{\nu}$. To this end, the following parameters are used: $C_{\nu} = 50$, $C_p = 1$, $C_z = 0.02$, $C_s = 1$ and $\eta = 0.01$. Initially, let the outer radius $R(0) = 2.71$, and inner radius $R_1(0) = 0.0151$ is obtained by solving \eqref{rel:2d3}. Then initial condition $\Sigma(0,x)$ and $W(0,x)$ are chosen of the form \eqref{Sigma-R-2d-lim} and \eqref{W-R-2d-lim} respectively, where $\Omega_1(t) = [- R_1(0), R_1(0)]$ and $\Omega_2 = [-R(0), -R_1(0)] \cup [R_1(0), R(0)]$. The solutions are gathered in Fig.~\ref{fig:radial1}. Here a good match is observed between the numerical solution and analytical formula. 

\begin{figure}[h!] 
\includegraphics[width=0.32\textwidth]{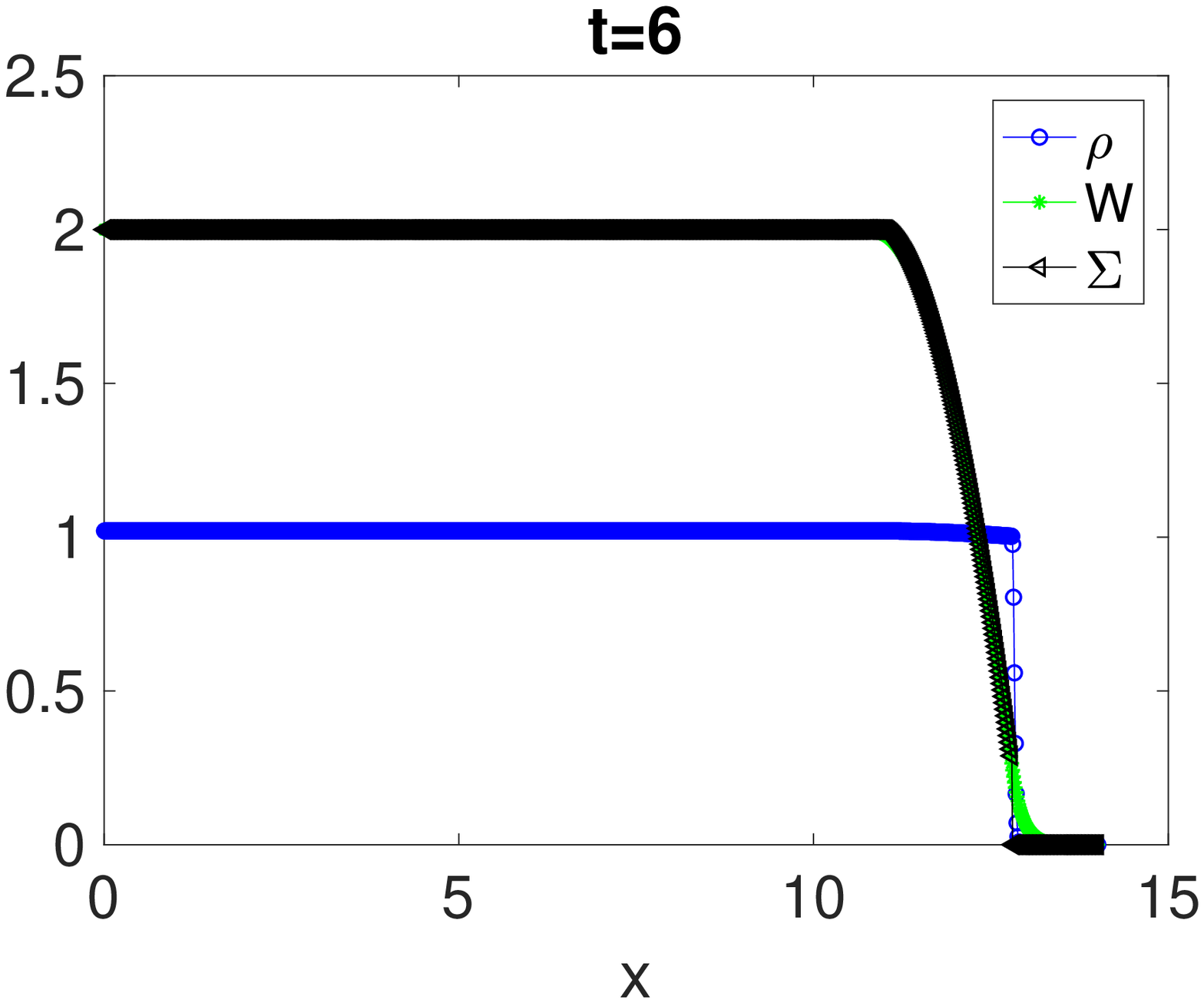}
\includegraphics[width=0.32\textwidth]{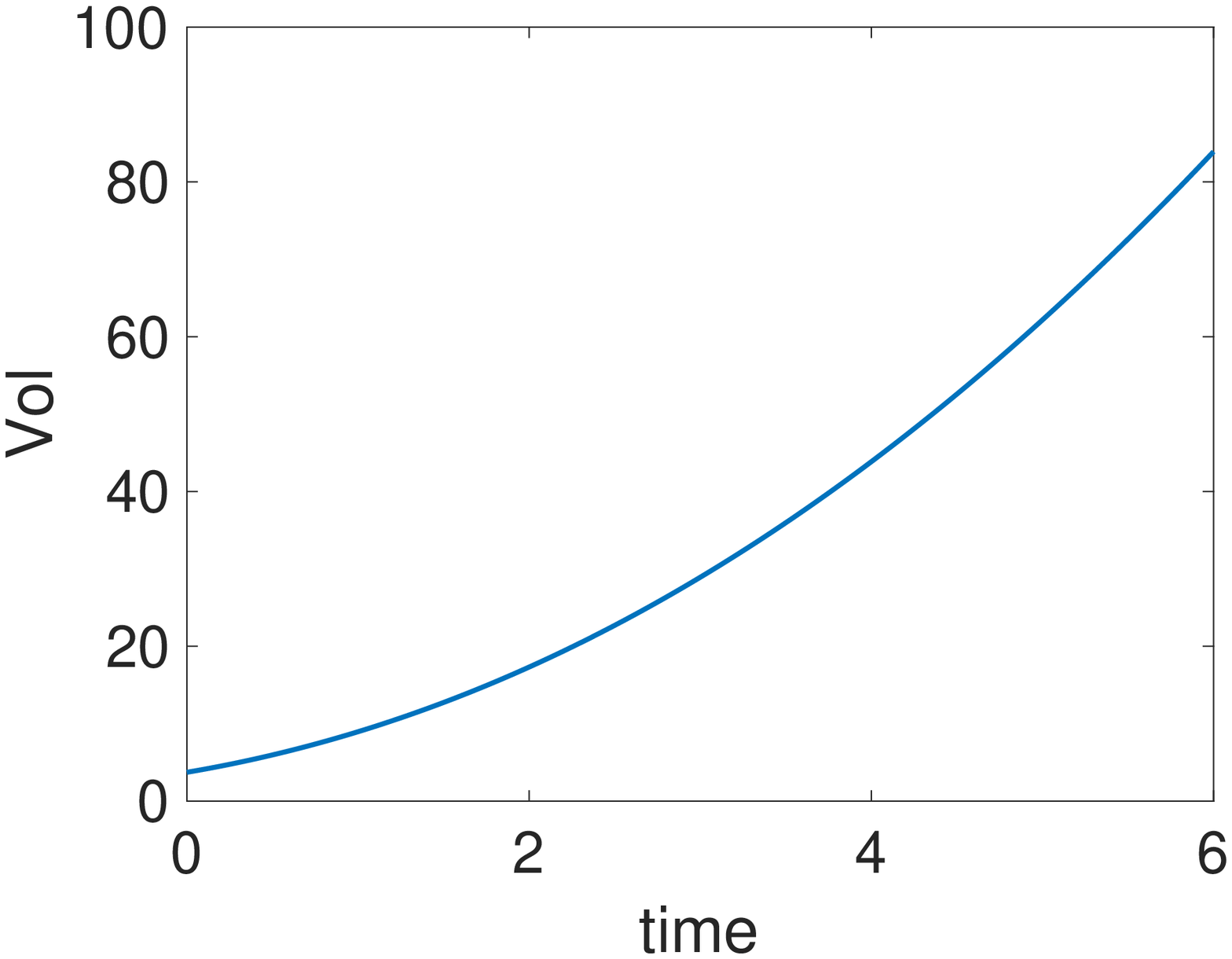}
\includegraphics[width=0.32\textwidth]{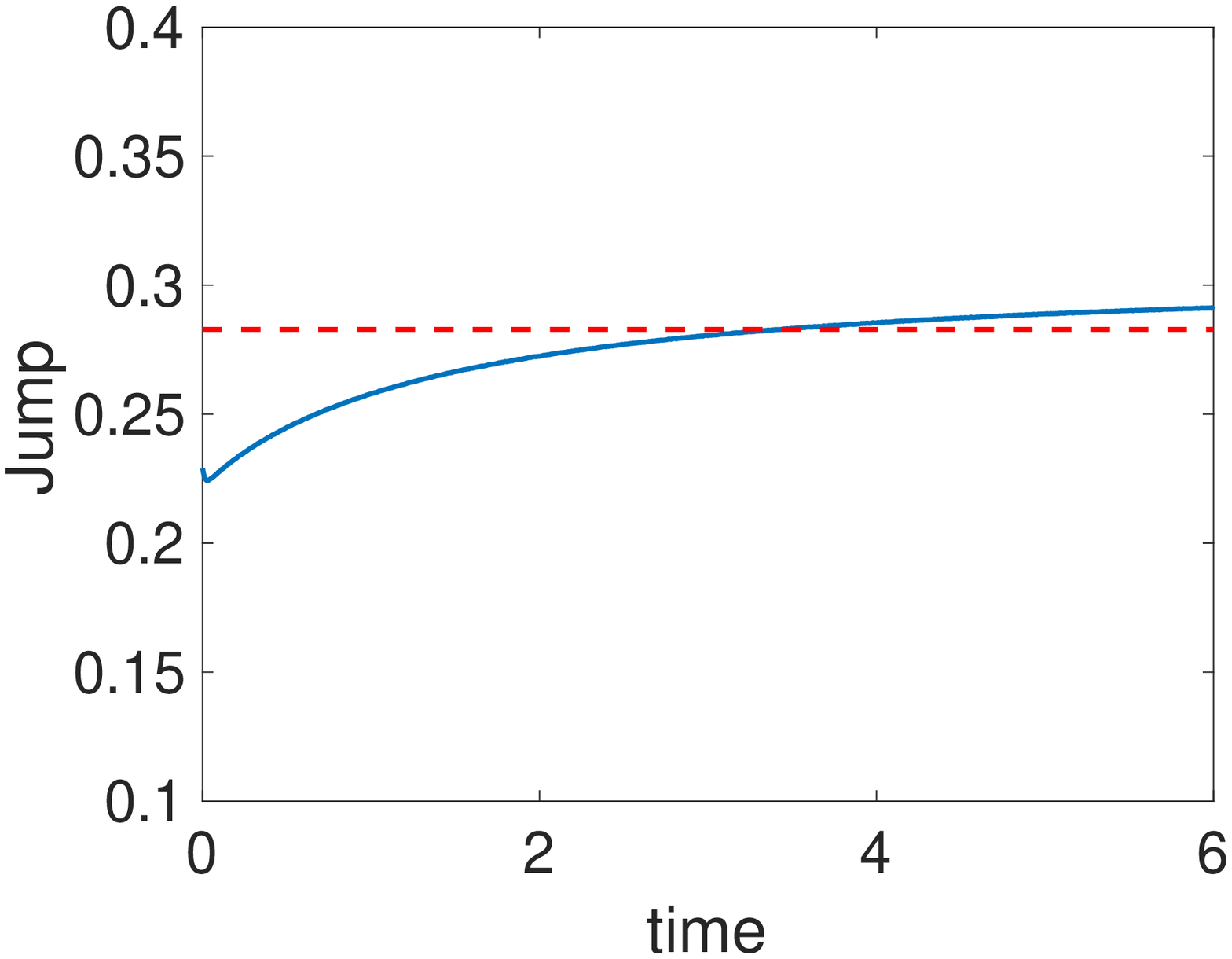}
\caption{2D radially symmetric case with initial data given by the limit model. Here $\Delta r = 0.0125$, $\Delta t = 1.067e-5$. The parameters used in this example are $C_{\nu} = 100$, $C_p = 2$, $C_z = 0.02$, $C_s = 1$ and $\eta = 0.0001$.}
\label{fig:radial12}
\end{figure}

\begin{figure}
\includegraphics[width=0.45\textwidth]{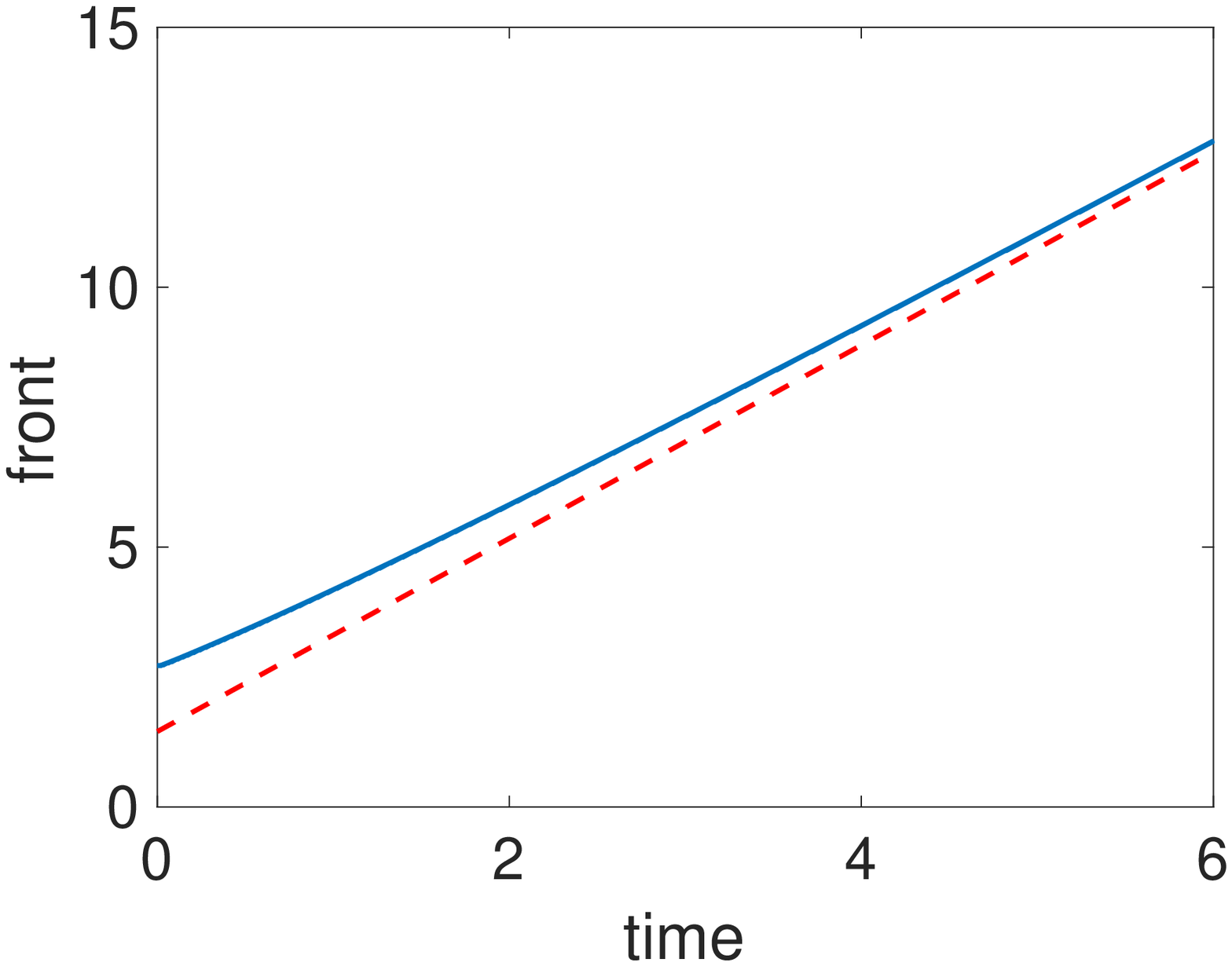}
\includegraphics[width=0.45\textwidth]{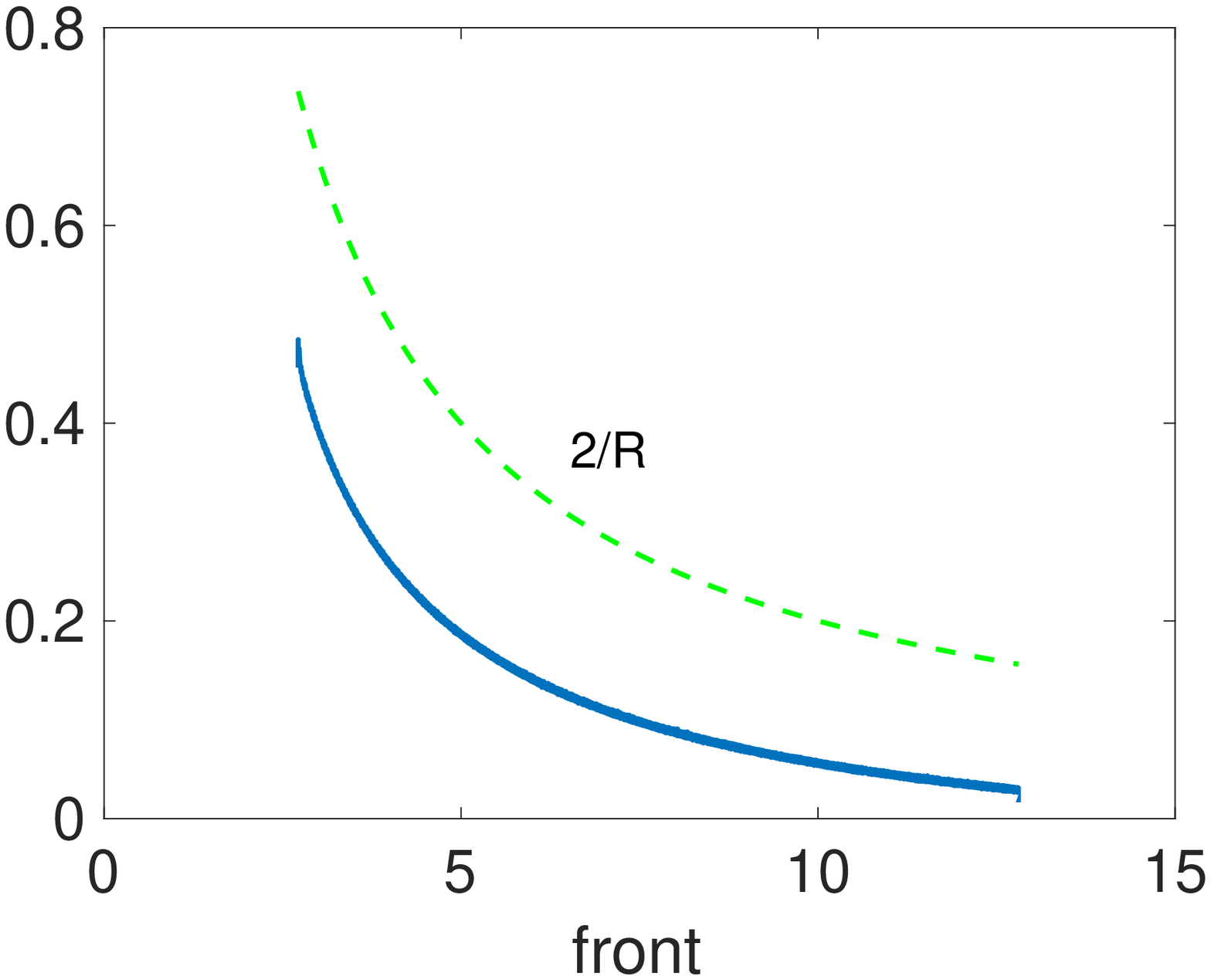}
\caption{Left: plot of the front versus time. The red dashed line represents slope $\sqrt{2C_p C_s} - \sqrt{C_z}$. Right: plot of $|\dot R(t)- (\sqrt{2C_p C_s} - \sqrt{C_z})$ versus front position, where $\dot R(t)$ denotes the front speed. Green dashed curve denotes $2/R$, where $R$ is the front position. Here $\Delta r = 0.0125$, $\Delta t = 1.067e-5$. The parameters used in this example are $C_{\nu} = 100$, $C_p = 2$, $C_z = 0.02$, $C_s = 1$ and $\eta = 0.0001$.}
\label{fig:radial13}
\end{figure}

Next, we check the jump in $\Sigma$, the volume of the tumor, and tumor invading front with respect to time in Figure \ref{fig:radial12}. The initial data is again chosen to be of the form \eqref{Sigma-1D} and \eqref{W-1D}  but with $R_1 = 1$, $R = 1.5$. The parameters are $C_{\nu} = 100$, $C_p = 2$, $C_z = 0.02$, $C_s = 1$ and $\eta = 0.0001$. We further check the convergence of propagation speed towards the limit in Figure \ref{fig:radial13}. Here the major difference compared to the 1D case is that we only observed the algebraic convergence, as denoted on the right of the plot. 

\section{Conclusion}

In this work, we explore the connections between a series of macroscopic models of tumor growth from the perspective of boundary propagation speeds. Prior to this work, only 1D traveling wave solutions are available, which yields a constant boundary moving speed. We give reassuring justification of the results in the traveling wave model since regardless of spatial dimension, the propagation speeds of radial symmetric solutions of the free boundary model all converge to that of the 1D traveling wave model. We also offer new observation that in multi-dimensional cases, the convergence of the propagation speed is algebraic and the curvature of the tumor profile, which is the reciprocal of the tumor radius, shows up in the first order correction to the boundary moving speed.  Between the cell density model and the free boundary model, we have numerically verified the the incompressible limit, which naturally implies the convergence of the propagation speed. But still, the rigorous convergence analysis is yet to be carried out, since the previous work only applies to the tumor growth models  coupled with the Darcy's law, but there should be no essential technical challenges when the Brinkman model is chosen. Besides, comprehensive numerical analysis and more general multi-dimensional implementation of the proposed numerical scheme are also worthy research topics. We shall pursue those directions in the future.  
\section*{Acknowledgments}
 J. Liu is partially supported by KI-Net NSF RNMS grant No. 11-07444 and NSF grant DMS 1812573. M. Tang is supported by Science Challenge Project No. TZZT2017-A3-HT003-F and NSFC 11871340. Z. Zhou is supported by NSFC grant No. 11801016. L. Wang is partially supported by NSF grant DMS-1903420 and NSF CAREER grant DMS-1846854.

\appendix

\section{Computing the regularized incompressible model in multi-dimensions}

\subsection{2D radial symmetric case}  \label{2d:cal}
Denote 
\[
\Omega_1^{\eta}= B_{R_1^\eta (t)}, \quad \Omega_1^{\eta} \cup \Omega_2^{\eta} = B_{R^\eta (t)},
\]
where $B_r$ denotes a ball centered at the origin with radius $r$, and assume $R^{\eta} (0)= R_0$. Hereafter, we shall first derive the relation between $R^{\eta}$ and $R_1^{\eta}$, and then the evolution equation for $R^{\eta}$. The derivation shares a lot in common with the one dimensional case, but results will have some subtle dependence on dimensions.

In $\Omega_1^\eta$, \eqref{limit000} becomes 
\[
-   \frac{C_S}{r} (rW_{r})_r = \frac{C_p -\Sigma}{\eta}, \quad 
-  \frac{C_z}{r} (rW_{r})_r  +W =\Sigma\,,
\]
which, by eliminating $\Sigma$, leads to
\[
-(\eta C_S+ C_z) \frac{1}{r} (rW_{r})_r +W=C_p.
\]
The symmetric assumption implies $W_r(0)=0$, and therefore the general solution of $W$ in $\Omega_1^\eta$ can be written as
\[
W(r)=C_p + A  I_0 \left(\frac{r}{\sqrt{\eta C_S+C_z}}\right), \quad  r\in \Omega_1^\eta\,,
\]
where $I_m(r)$ denotes the modified Bessel function of the first kind. Thus, the general solution of $\Sigma$ in $\Omega_1^\eta$ is given by
\[
\Sigma (r)=-  \frac{C_z}{r} (rW_{r})_r +W=C_p + \frac{A\eta C_S}{\eta C_S+ C_z} I_0 \left(\frac{r}{\sqrt{\eta C_S+C_z}}\right), \quad r \in \Omega_1^\eta.
\]
Note that at the boundary $r=R_1^{\eta}$ we have $\Sigma(R_1^{\eta})= C_p- \eta$, thus 
\begin{equation} \label{A-radial}
A= - \frac{\eta C_S+C_z}{C_S I_0 \left(\frac{R_1^{\eta}}{\sqrt{\eta C_S+C_z}}\right) }.
\end{equation}

In $\Omega_2^\eta$, \eqref{limit000} writes 
\[
-   \frac{C_S}{r} (rW_{r})_r = 1,\quad 
- \frac{C_z}{r} (rW_{r})_r +W =\Sigma, 
\]
which immediately leads to the general solution of $W$
\[
W(r)= - \frac{1}{4C_S} r^2 + a \ln r +b, \quad  r \in \Omega_2^\eta\,.
\]
By continuity of $W$ and $W_r$ at $r=R_1^{\eta}$, we get 
\begin{equation} \label{a-radial}
a= \frac{1}{2C_S}(R_1^\eta)^2 - R_1^{\eta} \frac{\sqrt{\eta C_S+C_z}}{C_S} \frac{ I_1 \left(\frac{R_1^{\eta}}{\sqrt{\eta C_S+C_z}}\right)}{ I_0 \left(\frac{R_1^{\eta}}{\sqrt{\eta C_S+C_z}}\right)},
\end{equation}
\begin{equation} \label{b-radial}
b=C_p - \eta - \frac{C_z}{C_S}+ \frac{(R^{\eta}_1)^2}{4C_S} -\frac{(R^{\eta}_1)^2 \ln R^\eta_1}{2C_S}  +R_1^{\eta} \ln R_1^\eta \frac{\sqrt{\eta C_S+C_z}}{C_S} \frac{ I_1 \left(\frac{R_1^{\eta}}{\sqrt{\eta C_S+C_z}}\right)}{ I_0 \left(\frac{R_1^{\eta}}{\sqrt{\eta C_S+C_z}}\right)}.
\end{equation}
And the solution of $\Sigma$ in $\Omega_2^\eta$ is given by
\[
\Sigma (r)=-  \frac{C_z}{r} (rW_{r})_r +W= - \frac{1}{4C_S} r^2 + a \ln r +b + \frac{C_z}{C_S}, \quad r \in \Omega_2 ^\eta \,.
\]

Finally, in $\Omega_3^\eta$, \eqref{limit000} reduces to 
\[
\Sigma=0, \quad 
-  \frac{C_z}{r} (rW_{r})_r +W =\Sigma.
\]
By assuming that $W$ decays at infinity, we have the following expression of $W$ in $\Omega_3^\eta$
\[
W(r)=d  K_0\left(\frac{r}{\sqrt{C_z}} \right), \quad r \in \Omega_3^\eta\,, 
\]
where $K_m(r)$ denotes the modified Bessel function of the second kind. The continuity of both $W$ $W_r$ at $R^\eta$ implies
\begin{equation} \label{d-radial}
d K_0\left(\frac{R^\eta}{\sqrt{C_z}} \right)=- \frac{1}{4C_S} (R^\eta)^2 + a \ln R^{\eta} +b.
\end{equation}
and
\begin{equation*}
 - \frac{d}{\sqrt{C_z}}K_1\left(\frac{R^\eta}{\sqrt{C_z}} \right) = - \frac{1}{2C_S} R^\eta + \frac{a}{R^\eta}.
\end{equation*}
In summary, the analytical representation of $\Sigma$ and $W$ are as follows
\begin{equation}\label{W-R-2d}
W(r) = \left\{ \begin{array}{cc} 
C_p + A  I_0 \left(\frac{r}{\sqrt{\eta C_S+C_z}}\right), & r\in \Omega_1^\eta
\\ - \frac{1}{4C_S} r^2 + a \ln r +b, &  r \in \Omega_2^\eta
\\ d  K_0\left(\frac{r}{\sqrt{C_z}} \right), & r \in \Omega_3^\eta
\end{array} \right.
\end{equation}
\begin{equation}\label{Sigma-R-2d}
\Sigma(r) = \left\{ \begin{array}{cc} 
C_p + \frac{A\eta C_S}{\eta C_S+ C_z} I_0 \left(\frac{r}{\sqrt{\eta C_S+C_z}}\right), & r\in \Omega_1^\eta
\\ - \frac{1}{4C_S} r^2 + a \ln r +b + \frac{C_z}{C_S}, &  r \in \Omega_2^\eta
\\ 0, & r \in \Omega_3^\eta\,,
\end{array} \right.
\end{equation}
where $A$, $a$, $b$ and $d$ are obtained from \eqref{A-radial}, \eqref{a-radial}, \eqref{b-radial} and \eqref{d-radial}, respectively.

\subsection{3D spherical symmetric case} \label{3d:cal}

For simplicity, we assume the problem is spherically symmetric in space, and we assume 
\[
\Omega_1^{\eta}= B_{R_1^\eta (t)}, \quad \Omega_1^{\eta} \cup \Omega_2^{\eta} = B_{R^\eta (t)},
\]
where $B_r$ denotes a ball centered at the origin with radius $r$.
And we assume the initial condition 
\[
R^{\eta} (0)= R_0.
\]
With the radial symmetric assumption, $W$ and $\Sigma$ are functions of only the radial variable $r$. The following calculations are similar to the 1D case, but we shall see some subtle effects of dimensions. 

First, we aim to derive the equation that link $R^{\eta}$ and $R_1^{\eta}$, and we plan to derive evolution equation that $R^{\eta}$ satisfies. 

In $\Omega_1^\eta$, the equations are
\[
-   \frac{C_S}{r^2} (r^2W_{r})_r = \frac{C_p -\Sigma}{\eta},
\]
\[
-  \frac{C_z}{r^2} (r^2W_{r})_r  +W =\Sigma.
\]
By eliminating $\Sigma$, we obtain
\[
-(\eta C_S+ C_z) \frac{1}{r^2} (r^2 W_{r})_r +W=C_p.
\]

The symmetric assumption implies $W'(0)=0$. Therefore, the general solution of $W$ in $\Omega_1^\eta$ is given by
\[
W=C_p + A \, i_0 \left(\frac{r}{\sqrt{\eta C_S+C_z}}\right),
\]
where $i_m(r)$ denotes the spherical modified Bessel function of the first kind. Thus, the general solution of $\Sigma$ in $\Omega_1^\eta$ is given by
\[
\Sigma=-  \frac{C_z}{r^2} (r^2 W_{r})_r +W=C_p + \frac{A\eta C_S}{\eta C_S+ C_z} i_0 \left(\frac{r}{\sqrt{\eta C_S+C_z}}\right).
\]

The boundary condition on at $r=R_1^{\eta}$ 
\[
\Sigma(R_1^{\eta})= C_p- \eta,
\]
leads to 
\begin{equation} \label{A-sph}
A= - \frac{\eta C_S+C_z}{C_S i_0 \left(\frac{R_1^{\eta}}{\sqrt{\eta C_S+C_z}}\right) }.
\end{equation}

In $\Omega_2^\eta$, the equations are
\[
-   \frac{C_S}{r^2} (r^2W_{r})_r = 1,
\]
\[
- \frac{C_z}{r^2} (r^2W_{r})_r +W =\Sigma.
\]
Obviously, the general solution of $W$ in $\Omega_2^\eta$ is given by
\[
W= - \frac{1}{6C_S} r^2 + a \frac 1 r +b.
\]
By continuity of $W$ and $W_r$ at $r=R_1^{\eta}$, we get 
\begin{equation} \label{a-sph}
a= -\frac{1}{3C_S}(R_1^\eta)^3 + (R_1^{\eta})^2 \frac{\sqrt{\eta C_S+C_z}}{C_S} \frac{ i_1 \left(\frac{R_1^{\eta}}{\sqrt{\eta C_S+C_z}}\right)}{ i_0 \left(\frac{R_1^{\eta}}{\sqrt{\eta C_S+C_z}}\right)},
\end{equation}
\begin{equation} \label{b-sph}
b=C_p - \eta - \frac{C_z}{C_S}+ \frac{(R^{\eta}_1)^2}{2C_S}   -R_1^{\eta}  \frac{\sqrt{\eta C_S+C_z}}{C_S} \frac{ i_1 \left(\frac{R_1^{\eta}}{\sqrt{\eta C_S+C_z}}\right)}{ i_0 \left(\frac{R_1^{\eta}}{\sqrt{\eta C_S+C_z}}\right)}.
\end{equation}
And the solution of $\Sigma$ in $\Omega_2^\eta$ is given by
\[
\Sigma=-  \frac{C_z}{r} (rW_{r})_r +W= - \frac{1}{6C_S} r^2 + a \frac 1 r +b + \frac{C_z}{C_S}.
\]

Finally, in $\Omega_3^\eta$, the equations are
\[
\Sigma=0,
\]
\[
-  \frac{C_z}{r} (rW_{r})_r +W =\Sigma.
\]
By assuming the decaying behavior at infinity, the general solution of $W$ in $\Omega_3^\eta$ is given by
\[
W=d  k_0\left(\frac{r}{\sqrt{C_z}} \right),
\]
where $k_m(r)$ denotes the modified Bessel function of the second kind. 

The continuity of $W$ at $R^\eta$ implies
\begin{equation}\label{d-sph}
d k_0\left(\frac{R^\eta}{\sqrt{C_z}} \right)=- \frac{1}{6C_S} (R^\eta)^2 + a \frac{1} {R^{\eta}} +b.
\end{equation}
And the continuity of $W_r$ at $R^\eta$ imposes a condition between $R_1^{\eta}$ and $R^\eta$
\begin{equation}\label{rel:3d2-2}
 - \frac{d}{\sqrt{C_z}}k_1 \left(\frac{R^\eta}{\sqrt{C_z}} \right) = - \frac{1}{3C_S} R^\eta - \frac{a}{(R^\eta)^2}. 
\end{equation}
In summary, the analytical representation of $\Sigma$ and $W$ are as follows
\begin{equation}\label{W-R-3d}
W(r) = \left\{ \begin{array}{cc} 
C_p + A \, i_0 \left(\frac{r}{\sqrt{\eta C_S+C_z}}\right), & r\in \Omega_1^\eta
\\ - \frac{1}{6C_S} r^2 + a \frac 1 r +b, &  r \in \Omega_2^\eta
\\d  k_0\left(\frac{r}{\sqrt{C_z}} \right), & r \in \Omega_3^\eta
\end{array} \right.
\end{equation}
\begin{equation}\label{Sigma-R-3d}
\Sigma(r) = \left\{ \begin{array}{cc} 
C_p + \frac{A\eta C_S}{\eta C_S+ C_z} i_0 \left(\frac{r}{\sqrt{\eta C_S+C_z}}\right), & r\in \Omega_1^\eta
\\ - \frac{1}{6C_S} r^2 + a \frac 1 r +b + \frac{C_z}{C_S}, &  r \in \Omega_2^\eta
\\ 0, & r \in \Omega_3^\eta\,,
\end{array} \right.
\end{equation}
where $A$, $a$, $b$ and $d$ are obtained from \eqref{A-sph}, \eqref{a-sph}, \eqref{b-sph} and \eqref{d-sph}, respectively.


\end{document}